\documentclass{article}
\usepackage{amsmath}
\usepackage{amsfonts}
\usepackage{amsthm}

\def\congruent{\equiv}

\newtheorem{Theorem}{Theorem}
\newtheorem{Lemma}{Lemma}
\newtheorem{Observation}{Observation}
\newtheorem*{Corollary}{Corollary to Lemma 2}
\newtheorem*{Corollary5}{Corollary to Lemma 5}
\newtheorem*{MainResult}{Main Result}
\newtheorem*{TheoremA}{Theorem A}

\begin{document}

\begin{center}
{\large\bf 
The generalized Pillai equation $\pm r a^x \pm s b^y = c$, II. 
}
\bigskip

Reese Scott

Robert Styer (correspondence author), Dept. of Mathematical Sciences, Villanova University, 800 Lancaster Avenue, Villanova, PA  19085--1699, phone 610--519--4845, fax 610--519--6928, robert.styer@villanova.edu
\end{center}


revised 20 Feb  2011 


\bigskip

\begin{abstract}  
We consider $N$, the number of solutions $(x,y,u,v)$ to the equation $ (-1)^u r a^x + (-1)^v s b^y = c $
in nonnegative integers $x, y$ and integers $u, v \in \{0,1\}$, for given integers $a>1$, $b>1$, $c>0$, $r>0$ and $s>0$.  When $(ra,sb)=1$, we show that $N \le 3$ except for a finite number of cases all of which satisfy $\max(a,b,r,s, x,y) < 2 \cdot 10^{15}$ for each solution; when $(a,b)>1$, we show that $N \le 3$ except for three infinite families of exceptional cases.  We find several different ways to generate an infinite number of infinite families of cases giving $N=3$ solutions.  
\end{abstract}   

\bigskip

\section{Introduction}
 
In this paper we consider $N$, the number of solutions $(x,y,u,v)$ 
to the equation  
$$ (-1)^u r a^x + (-1)^v s b^y = c \eqno{(1)}$$
in nonnegative integers $x, y$ and integers $u, v \in \{0,1\}$,
for given integers $a>1$, $b>1$, $c>0$, $r>0$ and $s>0$.  

In \cite{Sc-St3} we treated (1) with various additional restrictions on $x$, $y$, $u$, $v$, $a$, $b$.  In this paper, we treat (1) with no additional restrictions.  Brief histories of the problem are given in \cite{BL} and \cite{Sc-St3}, but see \cite{BBM} and \cite{W} for a much more extended history.  

Shorey \cite{Sh} showed that (1) has at most nine solutions in positive integers $(x,y)$ when $(u,v)=(0,1)$ and the terms on the left side of (1) are large relative to $c$.  

More recent results are given by (A), (B), (C), (D), and (E) which follow:

(A) $N \le 3$ when $x \ge 2$, $y \ge 2$, $(u,v)=(0,1)$, and $\gcd(ra,sb)=1$, except possibly when either $a$ or $b$ is less than $e^e$ (Le \cite{Le}).

(B) $N \le 2$ when $x \ge 2$, $y \ge 2$, $(u,v)=(0,1)$, and $\gcd(ra,sb)=1$, except possibly when $(a,b)$ is one of 23 listed pairs (Bo He and A. Togb\'e \cite{HT}).

(C) $N \le 2$ when $x \ge 1$, $y \ge 1$, $u$ and $v$ are unrestricted, and $\gcd(ra,sb)=1$, except for a finite number of cases which can be found in a finite number of steps and for which $\max(a,b,r,s,x,y) < 8 \cdot 10^{14}$ for  each solution \cite{Sc-St3}.  

(D) $N \le 3$ when $x \ge 1$, $y \ge 1$, $(u,v)=(0,1)$, and $\gcd(ra,sb)=1$, with no exceptions (Bo He and A. Togb\'e \cite{HT}).

(E)  $N \le 3$ when $x \ge 0$, $y \ge 0$, $(u,v) = (0,1)$, and $\gcd(ra,sb)$ is unrestricted, with no exceptions \cite{Sc-St3}.

In this paper we obtain several improvements on these results.  To state our main result, we will need a few preliminaries.  

We start with a few observations.

\begin{Observation}  
The choice of $x$ and $y$ uniquely determines the choice of $u$ and $v$.  
\end{Observation}

Following Observation 1, we will usually refer to a solution $(x,y)$.  

\begin{Observation}   
There are at most two solutions to (1) having the same value of $x$, similarly for $y$.  
\end{Observation}

In what follows we will often refer to a {\it set of solutions} to (1) which we will write as 
$$(a,b,c,r,s; x_1, y_1, x_2, y_2, \dots, x_N, y_N)$$
 and by which we mean the (unordered) set of ordered pairs $\{ (x_1, y_1), (x_2, y_2), \dots, (x_N, y_N) \}$, with $N > 2$, where each pair $(x_i,y_i)$ gives a solution to (1) for given integers $a$, $b$, $c$, $r$, and $s$.  

For Lemma 1 which follows, we define our use of the word {\it family}.  We say that two sets of solutions $(a,b,c,r,s;x_1, y_1, x_2, y_2, \dots, x_N, y_N)$ and $(A,B,C,R,S; X_1, Y_1, X_2, Y_2, \dots, X_N, Y_N)$ belong to the same family if $a$ and $A$ are both powers of the same integer, $b$ and $B$ are both powers of the same integer, there exists a positive rational number $k$ such that $kc=C$, and for every $i$ there exists a $j$ such that $k ra^{x_i} = R A^{X_j}$ and $ksb^{y_i} = S B^{Y_j}$, $1 \le i,j \le N$.  

\begin{Observation} 
It follows from the above definition of family that, if $(a,b,c,r,s; x_1, y_1, x_2, y_2, \dots, x_N, y_N)$ and $(A,B,C,R,S; X_1, Y_1, X_2, Y_2, \dots, X_N, Y_N)$ are in the same family with $kc = C$, then, by Observation 1, for every $i$ there exists a unique $j$ such that $k ra^{x_i} = R A^{X_j}$ and $ksb^{y_i} = S B^{Y_j}$, $1 \le i,j \le N$, and for every $j$ there exists a unique $i$ such that $k ra^{x_i} = R A^{X_j}$ and $ksb^{y_i} = S B^{Y_j}$, $1 \le i,j \le N$.    
\end{Observation}

\begin{Lemma} 
Every family contains a unique member $(a,b,c,r,s;x_1, y_1, x_2, y_2, \dots, x_N, y_N)$ with the following properties: $\gcd(r,s b)=\gcd(s, ra)=1$; $\min(x_1, x_2, \dots, x_N)=\min( y_1, y_2, \dots, y_N)=0$; and neither $a$ nor $b$ is a perfect power.  
\end{Lemma}

If a set of solutions has the properties listed in Lemma 1, we say it is in {\it basic form}.  

\begin{proof}[Proof of Lemma 1] 
Suppose a family contains a member $(a,b,c,r,s; x_1, y_1, x_2, y_2, \dots, x_N, y_N)$ in basic form, that is, with the properties $\gcd(r,s b)=\gcd(s, ra)=1$, $\min(x_1, x_2, \dots, x_N)=\min( y_1, y_2, \dots, y_N)=0$, and neither $a$ nor $b$ is a perfect power.  Then there must exist at least one $i$, $1 \le i \le N$, such that $(r a^{x_i}, s b^{y_i}) = 1$.  Assume $(A,B,C,R,S; X_1, Y_1, X_2, Y_2, \dots, X_N, Y_N)$ is another set of solutions in basic form belonging to the same family so that there exists a positive rational number $k$ such that $kc =C$.  Since for at least one $i$ we have no common factor dividing all of $r a^{x_i}$, $s b^{y_i}$, $c$, and since, by the definition of family, $k r a^{x_i}$, $k s b^{y_i}$, and $kc$, are all integers, we must have $k$ also an integer.  But then, if $k>1$, there does not exist a $j$ such that $(R A^{X_j}, S B^{Y_j})=1$ so that $(A,B,C,R,S; X_1, Y_1, X_2, Y_2, \dots, X_N, Y_N)$ cannot be in basic form.  So $k=1$, and there must exist a $j$ such that $r = R A^{X_j}$ and $X_j = \min(X_1, X_2, \dots, X_N)$; if $X_j \ne 0$, then again $(A,B,C,R,S; X_1, Y_1, X_2, Y_2, \dots, X_N, Y_N)$ cannot be in basic form; so $X_j=0$ and $r = R$.  In the same way, we can show that $k=1$ implies $s=S$.  Thus, we must have $a \ne A$ or $b \ne B$ since $(A,B,C,R,S; X_1, Y_1, X_2, Y_2, \dots, X_N, Y_N)$ is distinct from $(a,b,c,r,s; x_1, y_1, x_2, y_2, \dots, x_N, y_N)$, so that again $(A,B,C,R,S; X_1, Y_1, X_2, Y_2, \dots, X_N, Y_N)$ cannot be in basic form (by the last property listed in Lemma 1).  So we have shown that a given family contains at most one basic form.    

It remains to show that each family contains at least one basic form.  Now suppose a given family contains a set of solutions $(a,b,c,r,s; x_1, y_1, x_2, y_2, \dots, x_N, y_N)$ which is {\it not necessarily} in basic form.  
For any two solutions $(x_i, y_i)$ and $(x_j, y_j)$, $1 \le i,j \le N$, we have 
$$ r a^{\min{(x_i, x_j)}}(a^{|x_j - x_i|} + (-1)^\gamma ) = s b^{\min(y_i, y_j)} (b^{|y_j - y_i|} + (-1)^\delta )  \eqno{(2)} $$
where $\gamma, \delta \in \{0,1\}$.  
We can choose $(i,j)$, $ 1 \le i < j \le N$, so that $\min(x_i, x_j) = \min(x_1, x_2, \dots, x_N)= x_0$ and $\min(y_i, y_j) = \min(y_1, y_2, \dots, y_N)=y_0$.  For this choice of $(i,j)$, we have 
$$r a^{x_0}(a^t +(-1)^\gamma) = s b^{y_0}(b^w +(-1)^\delta) > 0 $$
where $t=|x_j - x_i|$, $w=|y_j - y_i|$, and $\gamma, \delta \in \{ 0, 1\}$.  
Let $g = \gcd(ra^{x_0}, s b^{y_0})$ and $h=\gcd(a^t +(-1)^\gamma, b^w +(-1)^\delta)$.  Then, taking $R= { r a^{x_0} \over g} = {b^w +(-1)^\delta \over h} $, $S= { s b^{y_0} \over g} = {a^t +(-1)^\gamma \over h}$, $C = {c \over g} $, we obtain a set of solutions $(a, b, C, R, S; x_1-x_0, y_1-y_0, \dots, x_N-x_0, y_N - y_0)$ in the same family as the original set.  In this set of solutions we can easily adjust $a$ and $b$, if necessary, so that neither $a$ nor $b$ is a perfect power.  The resulting set of solutions is in basic form except possibly when $\min(t, w)=0$, in which case without loss of generality we can take $t=0$, $w>0$, and $\gamma=0$ (note that since $c \ne 0$, we cannot have $t=w=0$).  We need to show $(a,S)=1$.  Assume $(a,S)>1$.  Then $S=2$, $2 | a$, and $R$ is odd.  For each $i$, $1 \le i  \le N$, we have $| R a^{x_i - x_0} \pm S b^{y_i - y_0}| = C$.  Choosing $i$ so that $x_i=x_0$, we get $C$ odd, while choosing $i$ so that $x_i > x_0$ (which we can do by Observation 2, noting that we have $N>2$ by the definition of a set of solutions), we get $C$ even, a contradiction which completes the proof of Lemma 1.  
\end{proof}

We are now ready to state the following

\begin{MainResult}
There are at most three solutions $(x,y,u,v)$ to (1) except for sets of solutions which are members of families which contain one of a finite number of basic forms which can be found in a finite number of steps.  

There are an infinite number of cases of three solutions to (1), even if we consider only sets of solutions in basic form.  
\end{MainResult}

This Main Result follows immediately from Theorems 1, 2, and 3, which are proved in Sections 3, 4, and 5.  In the case $(a,b)>1$, the present paper completely designates all exceptions.  In the case $(a,b)=1$, the present paper reduces the problem to a finite search; another paper \cite{Sc-St5} completes the search using not only the methods of \cite{GLS} and \cite{St1} as in previous work by the authors but also using LLL basis reduction and extensive computer algebra calculations.  The completion of the search in \cite{Sc-St5} proves Theorem A below, for which we need two more definitions. 

We define a { \it subset of a set of solutions} $(a,b,c,r,s; x_1, y_1, x_2, y_2, \dots, x_N, y_N)$ to be any set of solutions with the same $a$, $b$, $c$, $r$, $s$ and with all its pairs $(x,y)$ among the pairs $(x_i, y_i)$, $1 \le i \le N$.  Note that this subset may be (and, in our usage, usually is) the set of solutions $(a,b,c,r,s; x_1, y_1, x_2, y_2, \dots, x_N, y_N)$ itself.  

We define the {\it associate of a set of solutions } $(a,b,c,r,s; x_1, y_1, x_2, y_2, \dots, x_N, y_N)$ to be the set of solutions $(b,a,c,s,r; y_1, x_1, y_2, x_2, \dots, y_N, x_N)$.  

\begin{TheoremA} \cite{Sc-St5} 
Any set of solutions $(a,b,c,r,s; x_1, y_1, x_2, y_2, \dots, x_N, y_N)$ to (1) with $N > 3$ must be in the same family as a subset (or an associate of a subset) of one of the following: 
\setcounter{equation}{2}
\begin{align}
(3,2,1,1,2; 0,0,1,0,1,1,2,2) \notag \\
(3,2,5,1,2; 0,1,1,0,1,2,2,1,3,4) \notag \\
(3,2,7,1,2; 0,2,2,0,1,1,2,3)  \notag  \\
(5,2,3,1,2; 0,0,0,1,1,0,1,2,3,6)  \notag\\
(5,3,2,1,1; 0,0,0,1,1,1,2,3)   \\    
(7,2,5,3,2; 0,0,0,2,1,3,3,9)  \notag  \\
(6,2,8,1,7; 0,0,1,1,2,2,3,5) \notag  \\
(2,2,3,1,1; 0,1,0,2,1,0,2,0)  \notag  \\
(2,2,4,3,1; 0,0,1,1,2,3,2,4)  \notag
\end{align}
\end{TheoremA}

\section{Preliminary Lemmas}  

\begin{Lemma} 
If $(a,b,c,r,s; x_1, y_1, x_2, y_2, \dots, x_N, y_N)$ is a set of solutions to (1), then $\min(x_i, y_i) = 0$ for at most two choices of $i$, $1 \le i \le N$, except when $(a,b,c,r,s; x_1, y_1, x_2, y_2, \dots, x_N, y_N)$ is in the same family as a subset (or an associate of a subset) of one of the entries in the following list:
\begin{align*}
&(3,3,2,1,1;0,0,0,1,1,0), \\
&(5,2,3,1,2;0,0,0,1,1,0,1,2,3,6), \\
&(2,2,5,1,3;0,1,1,0,3,0), \\
&(2,2,3,1,1;0,1,0,2,1,0,2,0).
\end{align*}
\end{Lemma}

\begin{proof}
It suffices to show that, if $(a,b,c,r,s; x_1, y_1, x_2, y_2, \dots, x_N, y_N)$  is a set of solutions to (1) such that there are more than two values of $i$ for which $\min(x_i, y_i)=0$, the unique basic form in the same family is a subset (or an associate of a subset) of one of the listed exceptions in Lemma 2.  By Observation 3, we can assume this basic form or its associate has a subset $(A,B,C,R,S; X_1, Y_1, X_2, Y_2, X_3, Y_3)$ for which either 
$$X_1 = Y_1 = X_2 = Y_3 = 0, X_3>0, Y_2>0, \eqno{(4)}$$
or
$$X_1=Y_2=Y_3=0, X_3 > X_2 > 0, Y_1 > 0,  \eqno{(5)}$$
where $X_3 > X_2$ follows from Observation 1.  

If (4) holds, then, considering (2) with $(A,B,R,S, X_i, Y_i, X_j, Y_j)$ replacing $(a,b,r,s, x_i, y_i, x_j, y_j)$ and with $(i,j)=(1,2)$ and $(1,3)$, respectively, and noting $(R,S) = 1$ by the definition of basic form, we obtain $S \le 2$ and $R \le 2$ (note that neither side of (2) can  be zero).  Suppose $R=S=1$.  Then $C=2$ and, considering the solution $(X_2, Y_2)$, we get $B=3$ and ${Y_2}=1$, and, considering the solution $(X_3, Y_3)$, we get $A=3$ and $X_3=1$, giving the first exceptional case listed in Lemma 2; clearly there are no further solutions.  If $R \ne S$, then, by symmetry we can take $R=1$ and $S=2$.  Then $C=1$ or 3.  If $C=1$, then, considering the solution $(X_2, Y_2)$, we have a contradiction to (4).  So $C=3$, and, considering the solution $(X_2, Y_2)$, we get $B=2$ and $Y_2 = 1$, and, considering the solution $(X_3, Y_3)$, we get $A=5$ and $X_3 = 1$; for this choice of $(A,B,C,R,S)$, any further solutions $(X,Y)$ must have $\min(X,Y) > 0$, so that, by Theorem 1 of \cite{Sc-St2}, the only further solutions are $(X,Y) = (1,2)$ and $(3,6)$, giving the second exceptional case listed in the formulation of Lemma 2.  

Now suppose (5) holds.  Considering (2) as before with $(i,j)=(2,3)$, we get $R=1$, $A = 2$, and $X_2 = 1$.  Then considering (2) with $(i,j) = (1,2)$, we get either
$$S(B^{Y_1} + (-1)^\delta) = 3  \eqno{(6)}$$
or
$$S (B^{Y_1} - 1) = 1. \eqno{(7)}$$
Suppose (6) holds with $S=3$.  Then $B=2$, $Y_1 = 1$, and $\delta = 1$.  Considering the solution $(X_2, Y_2)$, we get $C=1$ or 5.  Considering the solution $(X_1, Y_1)$, we get $C = 7$ or 5.  So $C=5$, so that $X_3 = 3$, giving the third exceptional case listed in the formulation of Lemma 2; clearly no further solutions are possible.  

Now suppose (6) holds with $S=1$.  Then $(B, Y_1, \delta) = (2,1,0)$ or $(2, 2, 1)$.  $(B, Y_1, \delta) = (2,1,0)$ requires $C=1$, making the solution $(X_3, Y_3)$ impossible by Observation 1.  $(B, Y_1, \delta) = (2,2,1)$ requires $C=3$ and $X_3 = 2$.  In this case there is a fourth solution $(X,Y)=(0,1)$, giving the fourth exceptional case listed in the formulation of Lemma 2; clearly no further solutions are possible.  

If (7) holds, then $S=1$, $B=2$, and $Y_1=1$.  As in the immediately preceding case, we must have $C=3$ and $X_3 = 2$, and there exists the further solution $(X,Y)=(0,2)$, again giving the fourth exceptional case in the formulation of Lemma 2, which has no further solutions.  
\end{proof}

From Lemma 2 we immediately have the following: 

\begin{Corollary}
If a set of solutions to (1) is not in the same family as a subset (or an associate of a subset) of one of the entries listed in Lemma 2, and, further, if this set of solutions has at least one $x$ value equal to zero and at least one $y$ value equal to zero, then, letting $(a,b,c,r,s; x_1, y_1, x_2, y_2, \dots, x_N, y_N)$ be this set of solutions or its associate, we can assume one of the following holds: 
$$x_1 = y_1 = y_2 = 0, x_i > 0 {\rm \ for \ } i>1, y_i > 0 {\rm \ for \ } i > 2, \eqno{(8)}$$
$$x_1 = y_2 = 0, x_i > 0 {\rm \ for \ } i>1, y_1 > 0, y_i > 0 {\rm \ for \ } i > 2, \eqno{(9)}$$ 
$$x_1 = y_1 = 0, x_i > 0 {\rm \ for \ } i>1, y_1 > 0 {\rm \ for \ } i > 1. \eqno{(10)}$$ 
\end{Corollary}

\begin{Lemma}  
Let $(a,b,c,r,s; x_1, y_1, x_2, y_2, \dots, x_N, y_N)$ be a set of solutions to (1) for which $\gcd(a,b) > 1$ and $\min(x_1, x_2, \dots, x_N) = \min(y_1, y_2, \dots, y_N) = 0$.  If $(a,b,c,r,s; x_1, y_1, x_2, y_2, \dots, x_N, y_N)$ is not in the same family as a subset (or an associate of a subset) of one of the entries listed in Lemma 2, then we must have (10).  
\end{Lemma}

\begin{proof}
Let $(a,b,c,r,s; x_1, y_1, x_2, y_2, \dots, x_N, y_N)$ be a set of solutions to (1) for which $\gcd(a,b) > 1$ and $\min(x_1, x_2, \dots, x_N) = \min(y_1, y_2, \dots, y_N) = 0$.  Suppose $(a,b,c,r,s; x_1, y_1, x_2, y_2, \dots, x_N, y_N)$ is not in the same family as a subset (or an associate of a subset) of one of the entries listed in Lemma 2.  If Lemma 3 holds for a given set of solutions, it also holds for the associate of that set of solutions, so, by the Corollary to Lemma 2, we can assume that (8), (9), or (10) holds.  Suppose that (8) or (9) holds.  Then by Observation 3 we can assume $(a,b,c,r,s; x_1, y_1, x_2, y_2, \dots, x_N, y_N)$ is a member of a family containing a basic form which has as a subset a basic form $(A,B,C,R,S; X_1, Y_1, X_2, Y_2, X_3, Y_3)$ for which either 
$$X_1 = Y_1 = Y_2 = 0, X_2>0, X_3>0, Y_3 > 0  \eqno{(11)}$$
or 
$$X_1 = Y_2 = 0, X_2 > 0, X_3>0, Y_1 > 0, Y_3 > 0. \eqno{(12)}$$ 
Regardless of which of (11) or (12) holds, we have, considering (2) with $(i,j)=(2,3)$, 
$$R A^{\min(X_2, X_3)} ( A^{|X_3-X_2|} +(-1)^\gamma ) = S (B^{Y_3} +(-1)^\delta). \eqno{(13)}$$
Let $p$ be a prime dividing both $A$ and $B$.  Then $p$ divides the left side of (13) but not the right side of (13), since $(RA, S)=1$ by the definition of basic form.  This contradiction proves Lemma 3.
\end{proof}

\begin{Observation} 
If $p$ is a prime such that $p \mid a$ and $p \mid b$, and $t$ is the greatest integer such that $p^t \mid c$, then, if $(x,y)$ is a solution to (1), we must have $\min(x,y) \le t$.  
\end{Observation}

\begin{Lemma}  
Let $(a,b,c,r,s; x_1, y_1, x_2, y_2, \dots, x_N, y_N)$ be a set of solutions satisfying (10).  Then no two of $x_1$, $x_2$, \dots, $x_N$ are equal and no two of $y_1$, $y_2$, \dots, $y_N$ are equal, except when $(a,b,c,r,s; x_1, y_1, x_2, y_2, \dots, x_N, y_N)$ is in the same family as a subset (or an associate of a subset) of one of the following:
\begin{align*}
&(2,2,6,1,5; 0,0,2,1,4,1), \\
&(3,3,3,1,2; 0,0,1,1,2,1), \\
&(2,6,2,1,1; 0,0,2,1,3,1), \\
&(2,2,4,1,3; 0,0,3,2,4,2,1,1).
\end{align*}
\end{Lemma}

\begin{proof}
By Observation 3, any set of solutions contradicting Lemma 4 must occur in the same family as a basic form (or the associate of a basic form) which has as a subset a basic form $(A,B,C,R,S;\allowbreak  X_1, Y_1, X_2, Y_2, \allowbreak X_3, \allowbreak Y_3)$ for which 
$$X_1 = Y_1 = 0, X_3 > X_2>0 , Y_2 = Y_3 > 0.  \eqno{(14)}$$
To prove Lemma 4, it suffices to find each possible basic form $(A,B,C,R,S; X_1, Y_1, X_2, Y_2, \allowbreak X_3, \allowbreak Y_3)$ satisfying (14), determine its further solutions if any, and confirm that the list of exceptions in the formulation of Lemma 4 is complete.  

Consider (2) with $(A,B, R,S, X_i, Y_i, X_j, Y_j)$ replacing $(a,b, r,s, x_i, y_i, x_j, y_j)$.  Take $1 \le i,j \le 3$.  Combining (2) for $(i,j)=(1,2)$ with (2) for $(i,j)=(1,3)$, we get
$${ A^{X_3} + (-1)^{\gamma_3} \over A^{X_2} +(-1)^{\gamma_2} } = { B^{Y_2} +(-1)^{\delta_3} \over B^{Y_2} +(-1)^{\delta_2} }, \eqno{(15)}$$
where $\gamma_2$, $\gamma_3$, $\delta_2$, $\delta_3$ are in the set $\{0,1\}$.  Noting we cannot have $R A^{X_i} + S B^{Y_i} =C$ for $i>1$, we see that we must have  
$$\vert \gamma_2 - \delta_2 \vert = \vert \gamma_3 - \delta_3 \vert. \eqno{(16)}$$
The left side of (15) must be greater than or equal to 1, with equality when and only when $A=2$, $X_2 = 1$, $X_3 = 2$, $\gamma_2 = 0$, and $\gamma_3 = 1$.  Therefore, by (16), we have 
$$( \delta_2, \delta_3) = (1,0), \eqno{(17)}$$
so that $\gamma_2 \ne \gamma_3$ and, considering the right side of (15) as a reduced fraction and letting $m$ and $n$ be positive integers, we find that the right side of (15) must equal either ${ n+2 \over n}$ (when $B$ is even) or ${m +1 \over m}$ (when $B$ is odd).    If $A^{X_3 - X_2} \ge 5$, then both sides of (15) are greater than 3, which is impossible.  So we must have $2 \le A^{X_3-X_2} \le 4$.  

Assume first $A^{X_3-X_2} = 4$.  If $A^{X_2} \ge 8$, then both sides of (15) are greater than 3 which is impossible.  If $A^{X_2} = 2$, then the left side of (15) must be either 9 or $7/3$, neither of which is a possible value for the right side of (15).  So $A^{X_2} = 4$.  If $(\gamma_2, \gamma_3) = (1,0)$ then the left side of (15) equals $17/3$, which is again impossible.  So we are left with $A=2$, $X_2=2$, $X_3 = 4$, $\gamma_2 = 0$, $\gamma_3 = 1$, $B=2$, $Y_2=1$, $\delta_2 = 1$, and $\delta_3=0$.  Considering (2) with $(i,j)=(1,2)$ and recalling $(R,S)=1$ since we are dealing with a basic form, we find $R=1$, $S=5$, and $C=6$, giving the first exceptional case in the formulation of the lemma, which has no further solutions by Observation 4.  

Now assume $A^{X_3-X_2} = 3$.  If $A^{X_2} \ge 9$, then the left side of (15) is greater than 2 but not equal to 3, which is impossible by (17).  So we have $A=3$, $X_2=1$, and $X_3 = 2$.  Since neither side of (15) can be greater than 3, we must have $(\gamma_2, \gamma_3) = (0,1)$, giving $B=3$ and $Y_2=1$.  As in the previous paragraph, we find $R$, $S$, and $C$, obtaining the second exceptional case in the formulation of the lemma; applying Observation 4, we find that there are no further solutions.  

It remains to consider $A^{X_3-X_2} = 2$.  Consider first the case $(\gamma_2, \gamma_3) = (0,1)$.  If $A^{X_2} \ge 16$, the left side of (15) is greater than $5/3$ and less than 2, which is impossible by (17).  If $X_2=1$, the left side of (15) is equal to $1$, contradicting (17).  So we are left with either $X_2 = 2$ or $X_2=3$, so that, respectively, $B^{Y_2} = 6$ or $B^{Y_2} =4$.  In each of these two cases, we find $R$, $S$, and $C$ as in the preceding cases: when $X_2=2$, we obtain the third exceptional case in the formulation of the lemma, which has no further solutions by Observation 4; when $X_2=3$, we obtain a subset of the fourth exceptional case, which has no further solutions by Observation 4.    

Finally consider the case $A^{X_3-X_2} =2$ with $(\gamma_2, \gamma_3) = (1,0)$.  If $A^{X_2} \ge 8$, then the left side of (15) is greater than 2 and less than 3, contradicting (17).  If $X_2=1$, the left side of (15) equals 5, again contradicting (17).  So we are left with $X_2=2$, so that $B^{Y_2} = 2$, and, proceeding as in the preceding paragraphs, we obtain $(A,B,C,R,S ; X_1, Y_1, X_2, Y_2, X_3, Y_3) = (2,2,2,1,3;0,0,2,1,3,1)$, which has no further solutions by Observation 4 and which is in the same family as $(2,2,4,1,3;1,1,3,2,4,2)$, which is a subset of the fourth exceptional case listed in the formulation of the lemma.      
\end{proof}

To obtain results similar to that of Lemma 4 for the cases (8) and (9), we will need the following

\begin{Lemma}  
Suppose $(ra,sb)=1$ and (1) has two solutions $(x_1, y_1)$ and $(x_2, y_2)$ with $x_1 = 0$ and $x_2 > 0$.  Then, if $a$ is even, $r$ must be even.  
\end{Lemma}

\begin{proof}
If $a$ is even, then the solution $(x_2,y_2)$ requires $sb$ odd and $c$ odd, so that the solution $(x_1,y_1)$ requires $r$ even.
\end{proof}

\begin{Corollary5}
Suppose $(ra,sb)=1$ and (1) has three solutions $(x_1, y_1)$, $(x_2, y_2)$, $(x_3, y_3)$ with $x_1 < x_2$ and $x_1 < x_3$.  Then $y_2 \ne y_3$.  
\end{Corollary5}

\begin{proof}
Suppose $(ra,sb)=1$ and (1) has three solutions $(x_1, y_1)$, $(x_2, y_2)$, $(x_3, y_3)$ with $x_1< x_2$ and $x_1 < x_3$.  Suppose $y_2 = y_3$.  Then considering (2) with 
$(i,j)=(2,3)$, we find $ra^{\min(x_2, x_3)} | 2$, so that $a=2$, $r=1$, $\min(x_2, x_3) = 1$, and $x_1=0$, contradicting Lemma 5.
\end{proof}

We are now ready to prove Lemmas 6 and 7 which follow.

\begin{Lemma} 
Let $(a,b,c,r,s; x_1, y_1, x_2, y_2, \dots, x_N, y_N)$ be a set of solutions for which (8) holds.  Then no two of $x_1$, $x_2$, \dots, $x_N$ are equal and, except for $y_1=y_2$, no two of $y_1$, $y_2$, \dots, $y_N$ are equal, unless $(a,b,c,r,s; x_1, y_1, x_2, y_2, \dots, x_N, y_N)$ is in the same family as a subset of one of the following:
\begin{align*}
&(3,5,2,1,1; 0,0,1,0,1,1,3,2), \\
&(3,2,1,1,2; 0,0,1,0,1,1,2,2), \\
&(5,2,3,1,2; 0,0,1,0,1,2,3,6),  \\
&(2,7,3,2,1; 0,0,1,0,1,1).  
\end{align*}
\end{Lemma}

\begin{proof}
Let $(a,b,c,r,s; x_1, y_1, x_2, y_2, \dots, x_N, y_N)$ be a set of solutions satisfying (8), and assume there exists at least one pair $(i,j)$ where $2 \le i < j \le N$ such that either $x_i = x_j$ or $y_i = y_j$.  Then by Observation 3 $(a,b,c,r,s; x_1, y_1, x_2, y_2, \dots, x_N, y_N)$ must be in the same family as a basic form $(A,B,C,R,S; X_1, Y_1, X_2, Y_2, \allowbreak  \dots, \allowbreak X_N, Y_N)$ for which $X_1 = Y_1=Y_2=0$, $X_i>0$ for $i>1$, $Y_i>0$ for $i>2$, and there exists at least one pair $(i,j)$, where $2 \le i<j \le N$, for which $X_i =X_j$ or $Y_i = Y_j$.  By Lemma 3 we can assume $(a,b) = (A,B)=1$.  By the definition of basic form, we have $(RA, SB)=1$.  By the Corollary to Lemma 5, we can assume $X_i \ne X_j$ and $Y_i \ne Y_j$ for every pair $(i,j)$ such that $\min(i,j) \ge 3$.  Thus we can assume without loss of generality that $X_2 = X_3$, so that $(A,B,C,R,S; X_1, Y_1, X_2, Y_2, \dots, X_N, Y_N)$ has as a subset a basic form $(A,B,C,R,S; X_1, Y_1, X_2, Y_2, X_3, Y_3)$ for which 
$$X_1=Y_1=Y_2=0, X_2=X_3 >0, Y_3>0.  \eqno{(18)}$$
It suffices to find each possible basic form $(A,B,C,R,S;X_1, Y_1, X_2, Y_2, X_3, Y_3)$ satisfying (18), determine its further solutions if any, and thus verify the list of exceptions in the formulation of Lemma 6.  

Assume (18) holds, and consider (2) with $(A,B,R,S,X_i, Y_i, X_j, Y_j)$ replacing $(a,b,r,s, x_i, y_i, x_j, y_j)$.  Take $1 \le i, j \le 3$.  Then considering (2) with $(i,j)=(2,3)$, we have $S \le 2$.  And considering (2) with $(i,j)=(1,2)$, we have $R \le 2$.  

If $R=S=1$, then, considering the solution $(X_1, Y_1)$ to (1), we find $C=2$.  Considering the solution $(X_2, Y_2)$, we find $A=3$ and $X_2=1$.  Considering the solution $(X_3, Y_3)$, we find $B=5$ and $Y_3=1$.  Clearly, when $(A,B,C,R,S)=(3,5,2,1,1)$, any further solution to 
$$(-1)^U R A^{X} +(-1)^V S B^{Y} = C, \eqno{(19)}$$
where $U,V \in \{0,1\}$, must have 
$$X>0, Y>0, \eqno{(20)}$$
and 
$$(U,V) \ne (0,0).  \eqno{(21)}$$
By Theorem 7 of \cite{Sc-St2}, the only further solution is $(X_4, Y_4)=(3,2)$, giving the first exceptional case listed in the formulation of Lemma 6.  

If $R=1$ and $S=2$, then, considering (19) with $(X,Y)=(X_1,Y_1)$, we have either $C=1$ or $C=3$.  If $C=1$ then, considering the solution $(X_2, Y_2)$, we find $A=3$ and $X_2 = 1$.  Considering the solution $(X_3, Y_3)$, we find $B=2$ and $Y_3=1$.  Any further solution must satisfy (20) and (21).  By Theorem 7 of \cite{Sc-St2}, the only remaining solution is $(X_4, Y_4)=(2,2)$, giving the second exceptional case in the formulation of Lemma 6.  
So we now consider $C=3$.  Considering the solution $(X_2,Y_2)$, we find $A=5$ and $X_2=1$. Considering $(X_3,Y_3)$, we find $B=2$ and $Y_3=2$.  Any further solution must satisfy (21), and, since we are dealing with (8), it suffices to consider only further solutions which satisfy (20).  Theorem 7 of \cite{Sc-St2} shows that the only further solution with $\min(X,Y)>0$ is $(X_4, Y_4)=(3,6)$, giving the third exceptional case in the formulation of Lemma 6.   

Finally, considering $R=2$ and $S=1$, we find $C=1$ or $C=3$ again.  If $C=1$, then, considering the solution $(X_2, Y_2)$, we find $X_2=0$, a contradiction.  So we must have $C=3$.  Considering the solution $(X_2, Y_2)$, we find $A=2$ and $X_2=1$.  Considering $(X_3, Y_3)$, we find $B=7$ and $Y_3=1$.  Any further solution must satisfy (20) and (21).  Once again we can use Theorem 7 of \cite{Sc-St2}, this time finding there are no further solutions; we obtain the fourth exceptional case listed in the formulation of Lemma 6.  
\end{proof}

\begin{Lemma}  
Let $(a,b,c,r,s; x_1, y_1, x_2, y_2, \dots, x_N, y_N)$ be a set of solutions for which (9) holds.  Then no two of $x_1$, $x_2$, \dots, $x_N$ are equal and no two of $y_1$, $y_2$, \dots, $y_N$ are equal, unless $(a,b,c,r,s; x_1, y_1, x_2, y_2, \dots, x_N, y_N)$ is in the same family as a subset (or an associate of a subset) of one of the following:
\setcounter{equation}{21}
\begin{align}
&(3,2,5,1,2;0,1,1,0,1,2,2,1,3,4), \\
&(5,3,4,1,1;0,1,1,0,1,2), \\ 
&(3,2,7,1,2;0,2,2,0,2,3,1,1), \\
&(5,2,3,1,2; 0,1,1,0,1,2,3,6), \\
&(2^{g} + (-1)^\epsilon, 2, 2^{g} - (-1)^\epsilon, 1, 2; 0,g-1,1,0,1,g).  
\end{align}
where $\epsilon \in \{0,1\}$ and $g >1$ is a positive integer.  
\end{Lemma}

\begin{proof}
Proceeding as in the proof of Lemma 6, we see that it suffices to find each possible basic form $(A,B,C,R,S; X_1, Y_1, X_2, Y_2, X_3, Y_3)$ with $(RA, SB)=1$ for which 
$$X_1=Y_2=0, X_2=X_3 > 0, Y_1>0, Y_3>0, \eqno{(27)}$$ 
determine further solutions if any, and thus verify the list of exceptions in the formulation of Lemma 7.  Note that we can use (27) without loss of generality, since, if Lemma 7 holds for the associate of a set of solutions, it holds for the set of solutions itself.       

Suppose (27) holds.  
Considering (2) with $(i,j)=(2,3)$, we get 
$$ R A^{X_2} = { S \over 2} \left( B^{Y_3} + (-1)^{\delta} \right) \eqno{(28)}$$
where $\delta \in \{0,1\}$, so that 
$$C = R A^{X_2} - (-1)^{\delta} S = S B^{Y_3} - R A^{X_2}. \eqno{(29)}$$
Note that (28) implies $S \le 2$.  

Combining (28) and (29), we find
$$C = {S \over 2} \left( B^{Y_3} - (-1)^\delta \right). \eqno{(30)}$$
Suppose $Y_1 \ge Y_3$.  Then $S B^{Y_1} \ge S B^{Y_3} > R A^{X_2} > R$, so that 
$C= S B^{Y_1} \pm R  > S B^{Y_3} - R A^{X_2} = C$, a contradiction.  So we have a positive integer $d = {Y_3-Y_1}$.  Considering the solution $(X_1, Y_1)$ and using (30), we find that one of the following three equations must hold:
$$R = { S \over 2} \left( \left( B^d - 2 \right) B^{Y_1} - (-1)^{\delta}\right), \eqno{(31a)}$$
$$R = { S \over 2} \left( \left( B^d + 2 \right) B^{Y_1} - (-1)^{\delta}\right), \eqno{(31b)}$$
$$R = { S \over 2} \left( \left( 2 - B^d \right) B^{Y_1} + (-1)^{\delta}\right). \eqno{(31c)}$$

Suppose (31a) holds.  Then, since ${ 2 \over S} R $ divides ${2 \over S} R A^{X_2} - {2 \over S} R$, we have, using (28) and (31a), 
$$ \left( B^d - 2 \right) B^{Y_1} - (-1)^\delta  \mid 2 B^{Y_1} + (-1)^\delta 2 > 0.  \eqno{(32)} $$
From this we get 
$$\left( B^d - 4 \right) B^{Y_1} \le 3, $$
so that $B^{d} \le 5$.  If $B^d = 5$, then $B^{Y_1} = 2$ or 3, which is impossible.  So we have $ 2 \le B^d \le 4$.  

If $B^d=2$, then, by Lemma 5, $S=2$, and by (31a), $R=1$ and $\delta = 1$; further, by (28) we have $A^{X_2} = 2^{Y_1 + 1} -1$ and by (30) we have $C= 2^{Y_1 + 1} + 1$.  It is a well known elementary result that we must have $X_2=1$.  Thus we find that $(A,B,C,R,S; X_1, Y_1, X_2, Y_2, X_3, Y_3)$ must be (26) with $\epsilon=1$ and $g=Y_1+1$.  To see if there are any further solutions, we consider solutions $(X_i, Y_i)$ with $i>1$ and apply Theorem 1 of \cite{Sc-St2}, noting that $c$ must be a Fermat number greater than 3; we find that the only cases with further solutions are given by subsets of (22).  

Now suppose $B^d=3$.  
Then, from (32), we get 
$$  3^{Y_1} - (-1)^\delta  \mid 2 \cdot 3^{Y_1} + (-1)^\delta 2 = 2 \cdot 3^{Y_1} - (-1)^\delta 2 + (-1)^\delta 4,   $$
so that $3^{Y_1} - (-1)^\delta \mid 4$, so that $Y_1=1$.  
If $\delta = 0$, we find, using (31a), (30), and (28), that $(A,B,C,R,S; \allowbreak X_1, Y_1, \allowbreak X_2, Y_2, \allowbreak X_3, Y_3)$ must be (23), which has no further solutions by Theorem 7 of \cite{Sc-St2} (noting that we cannot have $R A^{X_4} + S B^{Y_4} = C$).  If $\delta = 1$, we find that $(A,B,C,R,S; X_1, Y_1, X_2, Y_2, X_3, Y_3) = (2,3,5,2,1;0,1,1,0,1,2)$; by Theorem 1 of \cite{Sc-St2}, we see that there are only two further solutions, giving (22) with the roles of $A$ and $B$ reversed.  
   
Now suppose $B^d = 4$.  
From (32) we have 
$$  2 B^{Y_1} - (-1)^\delta  \mid 2 B^{Y_1} + (-1)^\delta 2 = 2 B^{Y_1} - (-1)^\delta + (-1)^\delta 3, $$
so that $2 B^{Y_1} - (-1)^\delta \mid 3$, so that, since $Y_1 > 0$,  we must have $(B, Y_1, \delta) = (2, 1, 0)$, from which we find that $(A,B,C,R,S; X_1, Y_1, X_2, Y_2, X_3, Y_3) =(3,2,7,3,2; 0,1,1,0,1,3)$, which has no further solutions by Theorem 1 of \cite{Sc-St2}, and which is in the same family as $(3,2,7,1,2;1,1,2,0,2,3)$, which is a subset of (24).  So it suffices to consider (31b) and (31c).     

Suppose (31b) holds.  Then $0 < {2 \over S} R A^{X_2} - {2 \over S} R = - 2 B^{Y_1} + (-1)^\delta 2 $, which is impossible.   

Finally, suppose (31c) holds.  Then $B^d = 2$, $S=2$, $R=1$, and $\delta = 0$.  Further, by (28) we have $A^{X_2} = 2^{Y_1+1} + 1 $ and by (30) we have $C = 2^{Y_1+1} - 1$.  Here we have the possibility $X_2=2$, which gives $Y_1=2$, $A=3$, and $C=7$ and we obtain $(A,B,C,R,S; X_1, Y_1, X_2, Y_2, X_3, Y_3) = (3,2,7,1,2; 0,2,2,0,2,3)$, which has as its only further solution $(X,Y)=(1,1)$ by Theorem 1 of \cite{Sc-St2}, so we obtain (24).  If $X_2 \ne 2$, then we have $X_2=1$ and find that $(A,B,C,R,S; X_1, Y_1, X_2, Y_2, X_3, Y_3)$ must be (26) with $\epsilon = 0$ and $g=Y_1+1$; we apply Theorem 1 of \cite{Sc-St2} as above to see that the only case of (26) with $\epsilon = 0$ and $g>1$ allowing further solutions is given by (25).  (Note that (25) also has the solution $(0,0)$ which does not need to be taken into account here, since we are dealing with (9).)  
\end{proof}

\begin{Lemma} 
Let $(a,b,c,r,s; x_1, y_1, x_2, y_2, \dots, x_N, y_N)$ be a basic form satisfying (8) with no two of $x_2$, $x_3$, \dots,  $x_N$ equal, and no two of $y_3$, \dots, $y_N$ equal.   Then $x_2 = \min(x_2, x_3, \dots, x_N)$ except when $(a,b,c, r,s; x_1, y_1, x_2, \allowbreak y_2, \allowbreak \dots, \allowbreak x_N, y_N )$ is a subset of either $(3,2,5,1,4; 0,0,2,0,1,1,3,3)$ or $(2,3,5,2,3; 0,0,2,0,1,1,4,2)$.  
\end{Lemma}

\begin{proof}
Let $(a,b,c,r,s; x_1, y_1, x_2, y_2, \dots, x_N, y_N)$ be a basic form satisfying (8) with no two of $x_2$, $x_3$, \dots,  $x_N$ equal, and no two of $y_3$, \dots, $y_N$ equal, and suppose $x_2 \ne \min(x_2, x_3, \dots, x_N)$.  Without loss of generality take $x_3 = \min(x_2, x_3, \dots, x_N)$.  Combining (2) for $(i,j)=(1,2)$ and (2) for $(i,j)=(1,3)$, we get 
$$ { a^{x_3} + (-1)^{\gamma_3} \over a^{x_2} +(-1)^{\gamma_2}} = { b^{y_3} + (-1)^{\delta_3} \over 2}, \eqno{(33)}$$
where $\gamma_2$, $\gamma_3$, and $\delta_3$ are in the set $\{0,1\}$.  
Since $0 < x_3 < x_2$, the left side of (33) is less than or equal to one, with equality when and only when $a=2$, $ x_2=2$, $x_3 = 1$, $\gamma_2=1$, and $\gamma_3 = 0$.  
The right side of (33) is an integer or half integer greater than or equal to $1/2$.  So we have both sides of (33) equal either to $1/2$ or to 1.  

If both sides of (33) equal $1/2$, we must have $b^{y_3} = 2$ and $\delta_3=1$.  If also $\gamma_3=1$, the left side of (33) is less than $1/2$, so we have $\gamma_3=0$.  If $a^{x_2 - x_3} =2$, the left side of (33) is greater than $1/2$; if $a^{x_2-x_3} \ge 4$, the left side of (33) is less than $1/2$.  So we have $a^{x_2-x_3} = 3$, in which case, if $x_3 >1$, the left side of (33) is less than $1/2$.  This leaves as the only possibility $a=3$, $x_3 = 1$, $x_2=2$, $\gamma_3 = 0$, $\gamma_2 = 1$, $b=2$, $y_3 = 1$, $\delta_3 = 1$.  Considering (2) with $(i,j)=(1,2)$, we find $r/s = 1/4$.  By the definition of basic form, $(r,s)=1$, so $r=1$, $s=4$, and $c=5$, giving $(a,b,c,r,s; x_1, y_1, x_2, y_2, x_3, y_3)=(3,2,5,1,4; 0,0,2,0,1,1)$; if there exists a further solution $(x_4, y_4)$ we must have $x_4>0$, so that Theorem 1 of \cite{Sc-St2}, in combination with the fact that $3+2=5$, shows that the only possible further solution is $(x_4,y_4)=(3,3)$, giving the first exceptional case listed in the formulation of the lemma.  

Now suppose both sides of (33) are equal to 1.  Then we must have $a=2$, $x_2=2$, $x_3=1$, $\gamma_2=1$, $\gamma_3 = 0$, $b=3$, $y_3=1$, $\delta_3 = 1$.  Considering (2) with $(i,j)=(1,2)$, we find $r/s = 2/3$, so that $r=2$, $s=3$, and $c=5$, giving $(a,b,c,r,s; x_1, y_1, x_2, y_2, x_3, y_3) = (2,3,5,2,3;0,0,2,0,1,1)$; by Theorem 1 of \cite{Sc-St2} there is only one possible further solution $(x_4, y_4)=(4,2)$, giving the second exceptional case listed in the lemma.  
\end{proof}

\begin{Lemma}  
Let $(a,b,c,r,s; x_1, y_1, x_2, y_2, \dots, x_N, y_N)$ be a basic form satisfying (9) with no two of $x_2$, $x_3$, \dots, $x_N$ equal and no two of $y_1$, $y_3$, $y_4$, \dots, $y_N$ equal.  Then $x_2 = \min(x_2, x_3, \dots, x_N)$ and $y_1 = \min(y_1, y_3, y_4, \dots, y_N)$, except when $(a,b,c,r,s; x_1, y_1, x_2, y_2, \dots, x_N, y_N)$ or its associate is $(3,2,7,1,2; 0,2,\allowbreak 2,0,\allowbreak 1,1)$.  
\end{Lemma}

\begin{proof}  Assume $(a,b,c,r,s; x_1, y_1, x_2, y_2, \dots, x_N, y_N)$ is basic form satisfying (9) with no two of $x_2$, $x_3$, \dots, $x_N$ equal and no two of $y_1$, $y_3$, $y_4$, \dots, $y_N$ equal, and assume further that $x_2 \ne \min(x_2, x_3, \dots, x_N)$.  Without loss of generality we can take $x_3 = \min(x_2, x_3, \dots, x_N)$.  Choose $i \in \{1,3\}$ so that $y_i = \max(y_1, y_3)$, and let $n = 1$ or 2 according as $i=3$ or 1.  Let $k = s/c$.  Then 
$$k c = s \ge r a^{x_2} - c \ge a^n r a^{x_i} - c \ge a^n (s b^{y_i} -c) -c \ge a^n (b^2 k c - c) -c = (a^n b^2 k - (a^n + 1)) c \eqno{(34)}$$
so that 
$$k \le {a^n +1 \over a^n b^2 - 1}. \eqno{(35)}$$
Note that 
$$c > \left({1+k \over a}\right) c = {c+s \over a} \ge {r a^{x_2} \over a} \ge r a^{x_3} \eqno{(36)}$$
 
Assume $i=3$.  Then
$$b \le { s b^{y_3} \over s b^{y_1} } \le {c + r a^{x_3}  \over c - r  } \le { c + { r a^{x_2} \over a} \over c - {r a^{x_2} \over a^2} } \le { c + { (1+k) c \over a} \over c - { (1+k) c \over a^2} } = { 1+{1+k \over a} \over 1 - {1+k \over a^2}}.  \eqno{(37a)}$$
If $a \ge 3$ then (35) gives $ k \le 4/11$ so that (37a) gives $b < 2$, a contradiction.  So $a=2$, (35) gives $k \le 3/17$, so that (37a) gives $b < 3$, again a contradiction.  

So we must have $i=1$, $n=2$, giving 
$$b \le {s b^{y_1} \over s b^{y_3} } \le { c + r  \over c - r a^{x_3} } \le { c + { r a^{x_2} \over a^2} \over c - { r a^{x_2} \over a} } \le  { c + { (1+k) c \over a^2} \over c - { (1+k) c \over a} } = {1+ {1+k \over a^2} \over 1 - {1+k \over a}}. \eqno{(37b)} $$
Suppose $a \ge 3$.  Then (35) gives $ k \le 2/7$ and (37b) gives $b \le 2$, with equality only when $a=3$ and $s=(2/7) c$, so that, since $(s,c)=1$, we have $s = 2$ and $c=7$, so that, considering the solution $(x_2, y_2)$ and noting $x_2 > 1$, we get $r=1$ and $x_2 = 2$, so that $x_3=1$, and, considering the solution $(x_1, y_1)$, we get $y_1=2$, so that $y_3 = 1$, giving the exceptional set of solutions in the formulation of the lemma.  

If $a=2$, (35) gives $k \le 1/7$ and (37b) gives $b \le 3$, with equality only when $s=(1/7)c$, so that $s=1$, $c=7$, and, considering the solution $(x_2, y_2)$ and noting that $x_2 >1$ and $ 2 \mid r$ (by Lemma 5), we get $r=2$ and $x_2 = 2$, so that $x_3=1$, and, considering the solution $(x_1, y_1)$, we get $y_1 = 2$, $y_3 = 1$, giving the set of solutions in the formulation of the lemma with the roles of $a$ and $b$ reversed.  

In each of these two exceptional cases, Theorem 1 of \cite{Sc-St2} shows there is only one possible further solution $(x,y)$, given by $16-9 = 7$, which violates the assumption that no two positive $x$ values are equal and no two positive $y$ values are equal.  So we must have $x_2 = \min(x_2, x_3, \dots, x_N)$ except when $(a,b,c,r,s; x_1, y_1, x_2, y_2, \dots, x_N, y_N) \allowbreak = (3,2,7,1,2; 0,2,2,0,1,1)$ or $(2,3,7,2,1; 0,2,2,0,1,1)$.   

The same argument shows $y_1 = \min(y_1, y_3, y_4, \dots, y_N)$ except when $(a,b,c,r,s; x_1, y_1, x_2, y_2, \dots, x_N, y_N) \allowbreak  = (2,3,7,2,1; 0,2,2,0,1,1)$ or $(3,2,7,1,2; 0,2,2,0,1,1)$. 
\end{proof}

\begin{Lemma}  
Let $(a,b,c,r,s; x_t, y_t, x_h, y_h, x_k, y_k)$ be a set of solutions for which $x_t < x_h < x_k$, $y_t < y_h$, $y_t < y_k$, and $y_h \ne y_k$.   
Then $y_h < y_k$, except when $(a,b,c,r,s; x_t, y_t, x_h, y_h, x_k, y_k)$ is in the same family as $(2,2,2,1,1; 0,0,1,2,2,1)$. 
\end{Lemma} 

\begin{proof} 
Assume $(a,b,c,r,s; x_t, y_t, x_h, y_h, x_k, y_k)$ is a set of solutions for which $x_t < x_h < x_k$, $y_t < y_h$, and $y_t < y_k$, and assume further $y_k < y_h$.  Combining (2) for $(i,j)=(t,h)$ with (2) for $(i,j)=(t,k)$, we have
$${ a^{x_k - x_t} +(-1)^{\gamma_2} \over a^{x_h-x_t} +(-1)^{\gamma_1} } = { b^{y_k-y_t} +(-1)^{\delta_2} \over b^{y_h-y_t} +(-1)^{\delta_1} },  \eqno{(38)}$$
where $\gamma_1, \gamma_2, \delta_1, \delta_2 \in \{0,1\}$.  
The left side of (38) must be greater than or equal to one, with equality when and only when $a=2$, $x_h-x_t=1$, $x_k-x_t=2$, $\gamma_1 = 0$, $\gamma_2=1$.  But the right side of (38) must be less than or equal to one, with equality when and only when $b=2$, $y_h-y_t=2$, $y_k-y_t=1$, $\delta_1=1$, $\delta_2=0$, so we have $(a,b,c,r,s; x_t, y_t, x_h, y_h, x_k, y_k)$ in the same family as the basic form $(2,2,2,1,1; 0,0,1,2,2,1)$.  
\end{proof}

\begin{Lemma} 
Let $(a,b,c,r,s; x_1, y_1, x_2, y_2, \dots, x_N, y_N)$ be a set of solutions for which $\min(x_1, \allowbreak  x_2, \allowbreak \dots, \allowbreak  x_N) \allowbreak  = \min(y_1, y_2, \allowbreak  \dots, \allowbreak  y_N)=0$, no two solutions $(x,y)$ have the same positive $x$, and no two solutions $(x,y)$ have the same positive $y$.  Let $i$ and $j$ be distinct integers, $1 \le i,j \le N$.  Then $0 < x_i < x_j$ implies $y_i < y_j$, and $0 < y_i < y_j$ implies $x_i < x_j$, except when $(a,b,c,r,s;x_1, y_1, x_2, \dots, x_N, y_N)$ is in the same family as a subset (or an associate of a subset) of one of the following:
\begin{align*}
&(3,2,7,1,2; 0,2,2,0,1,1), \\
&(3,2,5,1,2; 0,1,1,0,1,2,2,1,3,4), \\
&(2,2,2,1,1; 0,0,1,2,2,1), \\
&(2,2,5,1,3; 0,1,1,0,3,0), \\
&(2,2,3,1,1;0,1,0,2,1,0,2,0).
\end{align*} 
\end{Lemma}

\begin{proof}
Assume $(a,b,c,r,s; x_1, y_1, x_2, y_2, \dots, x_N, y_N)$ is a set of solutions for which $\min(x_1, x_2, \allowbreak  \dots, \allowbreak  x_N) = \allowbreak \min(y_1, y_2, \allowbreak  \dots, \allowbreak  y_N)=0$, no two solutions $(x,y)$ have the same positive $x$, and no two solutions $(x,y)$ have the same positive $y$.  

Note that, if Lemma 11 holds for a given set of solutions, it holds for the associate of that set of solutions.  Note also that, of the exceptional sets of solutions in Lemma 2, the first two cases (and any subset of the second case) have no members of their families contradicting Lemma 11, and the last two are listed in Lemma 11.  Thus, by the Corollary to Lemma 2, it suffices to show that Lemma 11 holds for $(a,b,c,r,s; x_1, y_1, x_2, y_2, \dots, x_N, y_N)$ under the assumption that (8), (9), or (10) holds.     

If (10) holds, let $t=1$, $2 \le h \le N$, $2 \le k \le N$, and take $x_h < x_k$.  Then we have (38) with the conditions given in Lemma 10, so that Lemma 10 applies to prove Lemma 11  unless $(a,b,c,r,s; x_1, y_1, x_2, y_2, \dots, \allowbreak x_N, y_N)$  has a subset in the same family as $(2,2,2,1,1; 0,0,1,2,2,1)$; since $x_t = y_t = 0$, this subset has no further solutions since $(2,2,2,1,1; 0,0,1,2,2,1)$ has no further solutions by Observation 4 (immediately preceding Lemma 4).  Thus $(a,b,c,r,s; x_1, y_1 , x_2, y_2,\dots, x_N, y_N)$ itself is in the same family as $(2,2,2,1,1; 0,0,1,2,2,1)$.  This proves Lemma 11 when (10) holds.      

Suppose (8) holds with $x_2 \ne \min(x_2, x_3, \dots, x_N)$.  Then, by Lemma 8 and Observation 3, $(a,b,c,r,s; \allowbreak x_1, y_1, \allowbreak x_2, y_2, \dots, \allowbreak x_N, y_N)$ must be in the same family as a subset of one of the basic forms $(3,2,5,1,4; 0,0, \allowbreak 2,0,\allowbreak 1,1,\allowbreak 3,3)$ or $(2,3,5,2,3; \allowbreak  0,0, \allowbreak 2,0,\allowbreak 1,1,\allowbreak 4,2)$, and therefore in the same family as a subset (or an associate of a subset) of the second exceptional case of the formulation of Lemma 11.  So we can assume $x_2 = \min(x_2, x_3, \dots, x_N)$.

Suppose (9) holds with either $x_2 \ne \min(x_2, x_3, \dots, x_N)$ or $y_1 \ne \min(y_1, y_3, \dots, y_N)$.  Then, by Lemma 9 and Observation 3, $(a,b,c,r,s;x_1, y_1, x_2, y_2, \dots, x_N, y_N)$ must be in the same family as the basic form  $(3,2,7,1,2;0,2,2,0,1,1)$ or its associate; since $(3,2,7,1,2; \allowbreak 0,2 \allowbreak ,2,0,1,1)$ is the first exception listed in Lemma 11, we can assume $x_2 = \min(x_2, x_3, \dots, x_N)$ and $y_1 = \min(y_1, y_3, \dots, y_N)$.  

Thus, when (8) or (9) holds it suffices to show that Lemma 11 holds for $ 3 \le i,j \le N$.  Take $t \le 2$, $ h \ge 3$, $k \ge 3$, and take $x_h < x_k$.  We have (38) as in Lemma 10, so that Lemma 10 proves Lemma 11 since, by Lemma 3, $\gcd(a,b)=1$ and the exceptional case of Lemma 10 is impossible.  
\end{proof}

\begin{Lemma}  
Assume a set of solutions is not in the same family as a subset (or an associate of a subset) of one of the following:
\setcounter{equation}{38}
\begin{align}
&(3,2,1,1,2; 0,0,1,0,1,1,2,2) \notag \\
&(5,3,2,1,1; 0,0,0,1,1,1,2,3) \notag \\
&(5,2,3,1,2; 0,0,0,1,1,0,1,2,3,6) \notag \\
&(3,2,5,1,2; 0,1,1,0,1,2,2,1,3,4) \notag \\
&(3,2,7,1,2; 0,2,2,0,1,1,2,3)  \notag \\
&(5,3,4,1,1; 0,1,1,0,1,2) \notag \\
&(7,2,3,1,2; 0,0,0,1,1,1)  \notag \\
&(2,2,2,1,1; 0,0,1,2,2,1)  \notag \\
&(2,2,4,3,1; 0,0,1,1,2,3,2,4)   \\    
&(6,2,2,1,1; 0,0,1,2,1,3)  \notag \\
&(2,2,6,5,1; 0,0,1,2,1,4) \notag \\
&(3,3,3,2,1; 0,0,1,1,1,2)  \notag \\
&(2,2,3,1,1; 0,1,0,2,1,0,2,0)  \notag \\
&(2,2,5,3,1; 0,1,1,0,0,3)  \notag \\
&(3,3,2,1,1; 0,0,0,1,1,0)  \notag \\
&(2^g+(-1)^\epsilon, 2, 2^g - (-1)^\epsilon, 1,2; 0,g-1, 1,0, 1,g).  \notag
\end{align}
where $\epsilon \in \{0,1\}$ and $g$ is an integer with $g+\epsilon > 3$.  Then, letting $(a,b,c,r,s; x_1, y_1, x_2, y_2, \dots, x_N, y_N)$ be this set of solutions or its associate, we can assume one of the following must hold:
$$x_1 < x_2 < \dots < x_N, y_1 = y_2 < y_3 < \dots < y_N, \eqno{(40)}$$
$$x_1 < x_2 < \dots < x_N, y_2 < y_1 < y_3 < \dots < y_N.  \eqno{(41)}$$
$$x_1 < x_2 < \dots < x_N, y_1 < y_2 < \dots < y_N, \eqno{(42)}$$
\end{Lemma}

\begin{proof}
Each exceptional case (and therefore each subset of each exceptional case) in the formulations of Lemmas 2, 4, 6, 7, and 11 is in the same family as a subset (or an associate of a subset) of one of the entries in (39).  On the other hand, the basic form $(A,B,C,R,S; X_1, Y_1, X_2, Y_2, \dots, X_N, Y_N)$ in the same family as $(a,b,c,r,s; x_1, y_1, x_2, y_2, \dots, x_N, y_N)$ must {\it not} be in the same family as a subset (or an associate of a subset) of one of the entries in (39).  Therefore, $(A,B,C,R,S; X_1, Y_1, X_2, Y_2, \dots, X_N, Y_N)$ cannot be in the same family as a subset (or an associate of a subset) of one of the exceptional cases of Lemmas 2, 4, 6, 7, or 11.  By Lemmas 2, 4, 6, 7, and 11, we can assume that we have one of the following:
$$0=X_1 < X_2 < \dots < X_N, 0=Y_1 = Y_2 < Y_3 < \dots < Y_N, \eqno{(43)}$$
$$0=X_1<X_2<\dots<X_N, 0=Y_2 <Y_1<Y_3<\dots <Y_N.  \eqno{(44)}$$
$$0 = X_1 < X_2 < \dots < X_N, 0 = Y_1 < Y_2 <\dots < Y_N, \eqno{(45)}$$
Applying Observation 3 proves the lemma.
\end{proof}

\section{The case $\gcd(a,b)>1$}  

\begin{Theorem}  
If $\gcd(a,b)>1$, then there are at most three solutions $(x,y,u,v)$ to (1), except for sets of solutions which are members of families containing one of the following basic forms (taking $r \le s$):
$(2,2,4,1,3;0,0,1,1,3,2,4,2)$, $(2,2,3,1,1;0,1,0,2,1,0,2,0)$, $(6,2,8,1,7;0,0,1,1,2,2,3,5)$.  

There are an infinite number of cases of three solutions to (1) when $\gcd(a,b)>1$, even if we consider only sets of solutions in basic form.  
\end{Theorem}

\begin{proof}
To handle the first paragraph of Theorem 1, it suffices to prove that there are no sets of solutions in basic form with $N>3$ when $\gcd(a,b)>1$, except for cases (or associates of cases) listed in the formulation of the theorem.  

Assume there exists a set of solutions $(a,b,c,r,s; x_1, y_1, x_2, y_2, \dots, x_N, y_N)$ in basic form with $N \ge 4$ and $\gcd(a,b)>1$.  The first two exceptional cases in the formulation of the lemma are the only entries (or associates of entries) in (39) with $\gcd(a,b)>1$ and $N > 3$, so from here on we can assume, by Lemma 12 and Lemma 3, that (42) holds, and we must have 
$$0=x_1 < x_2 <\dots < x_N, 0 = y_1 < y_2 < \dots < y_N.  \eqno{(46)}$$

Suppose for some prime $p$ we have $p^\alpha || a$ and $p^\beta || b$ with $\alpha, \beta >0$.  Recall that, by the definition of basic form, we have $(r,sb)=(s, ra) = 1$.  Let $i$ and $j$ be integers such that $1 \le i < j \le N$.  Then, since $x_i < x_j$ and $y_i < y_j$, we have from (2) 
$$ \alpha x_i = \beta y_i. \eqno{(47)}$$
Since $x_2>0$ and $y_2>0$, we have, for $i=2$ and 3,
$${x_i \over y_i} = {\beta \over \alpha}.$$
Let $\alpha_0 = \alpha/ \gcd(\alpha, \beta)$ and $\beta_0 = \beta/\gcd(\alpha, \beta)$.  Then for $i=2$ or 3 we have 
$$x_i = k_i \beta_0, y_i = k_i \alpha_0 \eqno{(48)}$$
where $k_i$ is a positive integer.  Now combining (2) for $(i,j)=(1,2)$ with (2) for $(i,j)=(1,3)$ we obtain 
$${ a^{\beta_0 k_3} + (-1)^{\gamma_3} \over a^{\beta_0 k_2} +(-1)^{\gamma_2} } = { b^{\alpha_0 k_3} +(-1)^{\delta_3} \over b^{\alpha_0 k_2} +(-1)^{\delta_2}} \eqno{(49)}$$
where $\gamma_2$, $\gamma_3$, $\delta_2$ and $\delta_3$ are in the set $\{0,1\}$.  

Since $c \le r+s$, we must have
$$u_i \ne v_i, 2 \le i \le 4, \eqno{(50)}$$ 
where $u_i$ and $v_i$ are the values of $u$ and $v$ in (1) when $(x,y)=(x_i, y_i)$.  

If $a^{\beta_0} = b^{\alpha_0}$, then we have $\gamma_3 = \delta_3$ and $\gamma_2 = \delta_2$ in (49) and, from (48), $a^{x_2} = b^{y_2}$.  Considering (2) with $(i,j)=(1,2)$ and letting $h=\gcd(a^{x_2} +(-1)^{\gamma_2}, b^{y_2} +(-1)^{\delta_2})$, we have $r = (b^{y_2} +(-1)^{\delta_2}) /h$ and $s=(a^{x_2} +(-1)^{\gamma_2})/h$ so that $r=s$.  But then, by (50), we have $c=0$, a contradiction.  So $a^{\beta_0} \ne b^{\alpha_0}$.  For convenience of notation let $A=a^{\beta_0}$, $B= b^{\alpha_0}$, $n=k_2$, $m=k_3$.  Then (49) becomes 
$${ A^m +(-1)^{\gamma_3} \over A^n +(-1)^{\gamma_2} } = {B^m +(-1)^{\delta_3} \over B^n + (-1)^{\delta_2}}. \eqno{(51)}$$
We have already shown $A=B$ implies $c=0$ so $A \ne B$.  Assume $A>B$.  (The remainder of the proof works also for $B>A$.  Indeed, since we have (46), the roles of $a$ and $b$ are interchangeable.)  

(51) implies 
$${A^m - 1 \over A^n +1} \le {B^m +1 \over B^n -1}. \eqno{(52)}$$
If $n \ge 2$ then (52) implies 
$$ {(B+2)^3 -1 \over (B+2)^2 +1} \le { B^3 + 1 \over B^2 -1} \eqno{(53)}$$
noting that $A \ge B+2$ since $\gcd(A,B)>1$.  
But (53) does not hold for $B \ge 2$, so we must have $n=1$ so that 
$$ {A^m -1 \over A+1}  \le {B^m +1 \over B-1}$$
which implies, when $m \ge 3$, 
$${(B+2)^3 -1 \over B+3} \le {B^3 + 1 \over B-1}.  \eqno{(54)}$$ 
But (54) does not hold for $B \ge 2$ so we must have $m=2$ and $n=1$.  By (50), we can let $D= |\gamma_3 - \delta_3| = |\gamma_2 - \delta_2|$ (note $D = 0$ when $u_1 \ne v_1$, and $D=1$ when $(u_1, v_1) = (0,0)$, where $u_1$ and $v_1$ are the values of $u$ and $v$ in (1) when $(x,y)= (x_1, y_1)$).  If $D=0$ then the only possible choice of signs is given by
$${ A^2 + 1 \over A -1} = {B^2+ 1 \over B-1}. \eqno{(55)}$$
Since the function $(w^2 + 1) / (w -1)$ is monotone increasing for $w > 1+ \sqrt{2}$, the only possible solution to (55) is $(A, B)= (3,2)$ which is not under consideration here since we are taking $\gcd(a,b)>1$.  

So $D=1$.  If $\gamma_2=\gamma_3$ then we have 
$$ {A^2+1 \over A+1} = {B^2 - 1 \over B-1} $$
which, since the right side is the integer $B+1$, requires $A=0$ or 1, a contradiction.  So we must have $\gamma_2 \ne \gamma_3$ so that
$$ {A^2 - 1 \over A+1} = {B^2+1 \over B-1} \eqno{(56)} $$
which, since the left hand side is the integer $A-1$, requires $B=2$ or 3.  If $B=2$ then from (56) we see that $A=6$ and we obtain $a=6$, $b=2$, $r=(2-1)/\gcd((2-1), (6+1))=1$, $s=(6+1)/\gcd((2-1), (6+1))=7$, $c=8$ (since $D=1$ implies $(u_1,v_1)=(0,0)$), $x_1 = 0$, $y_1 =0$, $x_2=1$, $y_2=1$, $x_3=2$, $y_3=2$.  We find there is a fourth solution $x_4=3$, $y_4=5$.  Noting that $2^3 || c$, we easily see there can be no further solutions (by Observation 4 immediately preceding Lemma 4).  We obtain the first of the exceptions in the formulation of Theorem 1.  If $B=3$, then, from (56), we find $A=6$ and we obtain $a=6$, $b=3$, $r= (3-1) / \gcd((6+1), (3-1)) = 2$, $s=(6+1)/\gcd((6+1), (3-1))= 7$, $c=9$, $x_1=0$, $y_1=0$, $x_2=1$, $y_2=1$, $x_3=2$, $y_3=2$.  Noting $3^2 || c$, we easily see there are no further solutions.  

This ends the proof of the first part of Theorem 1.  

For the proof of the second part of Theorem 1, we note that there are an infinite number of choices of $a$ no two of which have a common factor and for each of which there are an infinite number of choices of $m$ such that 
a set of solutions to (1) with $N=3$ is given by
$$(a,b,c,r,s; x_1, y_1, x_2, y_2, x_3, y_3) = (a, ta, {a (t+(-1)^{u+v+1}) \over h} , { ta + (-1)^{v} \over h}, { a +(-1)^{u} \over h}; 0,0,1,1,m+1,2) \eqno{(57)}$$
where $m \ge 0$ is an integer, $t = {a^m + (-1)^{v} \over a+(-1)^u}$ is an integer, $h= \gcd( ta + (-1)^v, a+(-1)^u)$, and $u$ and $v$ are in the set $\{ 0,1 \}$.  Thus there are an infinite number of families with $N = 3$.  (Note that we do not actually need the fact that there are an infinite number of choices of $m$ for each choice of $a$, although this fact is easily verified.)   
This completes the proof of Theorem 1.    
\end{proof}

A further example giving an infinite number of sets of solutions with $N=3$ is the following (closely related to (57)):
$$ (a,b,c,r,s; x_1, y_1, x_2, y_2, x_3, y_3) = (2, 4t, {4t+4 \over h_1} , { 4t+1 \over h_1}, { 3 \over h_1}; 0,0,2,1,m_1+2,2) \eqno{(58)} $$
where $m_1 \ge -1$ is an odd integer, $t = {2^{m_1} +1 \over 3}$, $h_1 = 3$ or 1 according as $m_1 \congruent 5 \bmod 6$ or not, and $u,v \in \{0,1\}$.

\section{The case $\gcd(a,b)=1$}  

Theorem 2, which follows, proves the first paragraph of the Main Result above for the case $(a,b)=1$ since any basic form for which $(a,b)=1$ must also satisfy $(ra, sb)=1$.

\begin{Theorem}  
When $(ra, sb)=1$, (1) has at most three solutions $(x,y, u,v)$ in nonnegative integers $x$, $y$, and $u,v \in \{0,1\}$, except for a finite number of sets of solutions $(a,b,c,r,s; x_1, y_1, x_2, y_2, \dots, x_N, y_N)$ all of which can be found in a finite number of steps.   

If (1) has more than three solutions $(x,y,u,v)$ when $(ra,sb)=1$, then
$$ \max(a,b,r,s,x,y) < 2 \times 10^{15} \eqno{(59)}$$
for each solution, and, further, 
$$ \max(r,s,x,y) < 8 \times 10^{14},  \min(\max(a,b), \max(b,c), \max(a,c)) < 8 \times 10^{14}  \eqno{(60)}$$
for each solution.  
\end{Theorem}

For the proof of Theorem 2, we require a few additional lemmas.

\begin{Lemma} 
Suppose $(ra,sb)=1$ and suppose (1) has four solutions $(x_1, y_1)$, $(x_2, y_2)$, $(x_3, y_3)$, $(x_4, y_4)$ with $x_1 < x_2 < x_3 < x_4$.  Let $Z = \max(x_4, y_1, y_2, y_3, y_4)$.   Then
$$a^{x_3-x_2} \le Z, s \le Z+1, $$
and, if $a>b$,
$$x_2 a^{x_3-x_2} \le Z. $$
\end{Lemma}

\begin{proof}
By the Corollary to Lemma 5, we see that no two of $y_2$, $y_3$, $y_4$ are equal.  Considering (2) with $(i,j)=(2,3)$ and $\delta= \delta_1 \in \{0,1\}$, we find $\left(b^{|y_3-y_2|} +(-1)^{\delta_1}\right)/\left(r a^{x_2}\right)$ is an integer prime to $a$.  Considering (2) with $(i,j)=(3,4)$ and $\delta=\delta_2 \in \{0,1\}$, we find $\left(b^{|y_4-y_3|} +(-1)^{\delta_2}\right)/\left(ra^{x_3}\right)$ is an integer.  There must be a least positive integer $n$ such that $(b^n \pm 1)/(r a^{x_2})$ is an integer prime to $a$.  Now we can apply Lemma 1 of \cite{Sc-St3} to get 
$$n { a^{x_3-x_2} \over 2^{g+h-1}} \le \max(y_3, y_4) \le Z,$$
where $g=1$ and $h=0$ unless $r$ is odd, $a \congruent 2 \bmod 4$, and $x_2=1$, in which case $x_1=0$.  But by Lemma 5 we cannot have $r$ odd and $a \congruent 2 \bmod 4$ when $x_1=0$, so that we can simply take $g+h-1=0$.  Also if $a>b$, then $n > x_2$.  

Finally, considering (2) with $(i,j)=(2,3)$ we find $s \le a^{x_3-x_2} + 1 \le Z+1$. 
\end{proof}

\begin{Lemma}  
Suppose $(ra,sb)=1$ and suppose (1) has 2 solutions $(x_1, y_1)$ and $(x_2, y_2)$, with $x_1 < x_2$.  Then, if $r>1$ or if $x_1 > 0$, 
$$r a^{x_2} > c/2,$$
and, if $r=1$ and $x_1 = 0$,
$$a^{x_2} > (c-2)/2.$$
\end{Lemma}

\begin{proof}
Suppose $ra^{x_2} \le c/2$.  Then $sb^{y_2} \ge c/2$ and also $ra^{x_1} < c/2$ so that $sb^{y_1} > c/2$.  But we have, considering (2) with $(i,j)=(1,2)$,  
$$r a^{x_1} ( a^{x_2-x_1} +(-1)^\gamma ) = s b^{\min(y_1, y_2)} \left( b^{|y_2-y_1|} + (-1)^\delta \right) > 0.$$
If $r>1$, or if $x_1>0$, this gives $c/2 \ge r a^{x_2} > a^{x_2-x_1} + 1 \ge s b^{\min(y_1, y_2)} \ge c/2$, a contradiction, so $r a^{x_2} > c/2$ when $r>1$ or $x_1>0$.  

If $r=1$ and $x_1=0$, suppose $a^{x_2} \le c/2 - 1$, so that $ s b^{y_2} \ge c/2 + 1$.  
Proceeding in the same way we obtain $c/2 \ge a^{x_2} + 1 \ge s b^{\min(y_1, y_2)} \ge c/2 + 1$, a contradiction.  So $a^{x_2} > c/2 - 1$ when $r=1$ and $x_1=0$.  
\end{proof}

\begin{Lemma}  
Suppose (1) has a solution $(x,y,u,v)$ for some $(a,b,c,r,s)$, and suppose we have the following conditions:
$$ \min(x,y) > 0, (u,v) \ne (0,0), (ra, sb) = 1. $$
Let $Z = \max(x,y)$, $J = \max(a,b)$, and $d = \min( r a^x, s b^y)$.  Then one of the following inequalities must hold:  
$$ Z < { \log(1+c/d) + \log(c) \over \log(2)} + 1.6901816335 \cdot 10^{10} \log(\max(r,s)) \log(J) \log(1.5 e Z), \eqno{(61)}$$
$$ Z < { \log(1+c/d) + \log(c) \over \log(2)} + 22.997 \left( \log\left({Z \over \log(2)}\right) + 2.405 \right)^2 \log(J). \eqno{(62)}$$
When $rs=1$, then either (62) holds or $Z < 2409.08 \log(J)$.  
\end{Lemma}

\begin{proof}
Suppose $(ra,sb)=1$.  Let $Z= \max(x,y)$ where $(x,y,u,v)$ is a solution in positive integers to (1).  Let $D=\max(r a^x, s b^y)$, $d = \min(r a^x, s b^y)$, $J = \max(a,b)$, and $j = \min(a,b)$.  Assume $(u,v) \ne (0,0)$  so that 
$$|r a^{x} - s b^{y}| = D - d = c.$$
Assume $rs>1$ and also assume there do not exist integers $a_0$, $b_0$, $m$, $n$, $t$, and $w$ such that $r = a_0^m$, $a=a_0^t$, $s=b_0^n$, and $b=b_0^w$.  Let 
$$\Lambda = |\log(r/s)+ x \log(a)-y \log(b)|=\log(D) - \log(d)= \log(1 + c/d)< c/d$$
so that 
$$\log(1/\Lambda) > \log(d) - \log(c). \eqno{(63)} $$

Now we can apply a result of Matveev \cite{Ma}, as given by Mignotte in Theorem 1 of \cite{Mi2}, with $(\alpha_1, \alpha_2, \alpha_3) = (r/s, a, b)$, $(A_1, A_2, A_3) =$ ($\log(\max(r,s))$, 
$\log(a)$, $\log(b)$), and $B = Z$ to get 
$$\log(1/\Lambda) < K A_1 A_2 A_3 \log(1.5 e B), \eqno{(64)}$$
 where $K = 1.6901816335 \cdot 10^{10}$.  Combining (63) and (64) we get 
$$
\log(d) < \log(c) + K  \log(\max(r,s)) \log(a) \log(b) \log(1.5 e Z). \eqno{(65)}$$
Also $\Lambda = \log(1+c/d) $, so that, adding $\Lambda$ to both sides of (65), we get
$$\log(D)< \log(1 + c/d) + \log(c)  + K  \log(\max(r,s)) \log(a) \log(b) \log(1.5 e Z).   $$
From this, noting that $Z \log(j) \le \log(D) $, we have
$$ Z \log(j) < \log(1 + c/d) + \log(c)+  K  \log(\max(r,s)) \log(a) \log(b) \log(1.5 e Z).$$
Dividing through by $\log(j)$ and noting $j \ge 2$, we get (61).

Now suppose $rs=1$.  Let $G = \max(x / \log(b), y / \log(a))$ and let $\Lambda = | x \log(a) - y \log(b) |$.  Using a theorem of Mignotte \cite{Mi1} as given in Section 3 of \cite{Be} and using the parameters chosen by Bennett in Section 6 of \cite{Be}, we see that we must have either 
$$G < 2409.08 \eqno{(66)}$$
or
$$\log(\Lambda) > -22.997 (\log(G) + 2.405 )^2 \log(a) \log(b). \eqno{(67)}$$
If (66) holds then $Z < 2409.08 \log(J)$, which implies (61).  
When (67) holds, proceeding as in the case $rs>1$ above and noting $G \le Z / \log(2)$, we obtain (62).

Suppose $rs>1$ and there exist integers $a_0$, $b_0$, $m$, $n$, $t$, and $w$ such that $r = a_0^m$, $a=a_0^t$, $s=b_0^n$, and $b=b_0^w$.  Then we can rewrite (1) as 
$$(-1)^u a_0^{tx+m} + (-1)^v b_0^{wy+n} = c.$$
From this we obtain, using the same method as above when $rs=1$, and letting $Z_0 = \max(tx+m, wy+n)$, and $J_0 = \max(a_0, b_0)$, 
$$Z \le Z_0 < 2409.08 \log(J_0) \le 2409.08 \log(J) \eqno{(68)}$$
or
$$Z_0 < { \log(1+c/d) + \log(c) \over \log(2)} + 22.997 \left( \log\left({Z_0 \over \log(2)}\right) + 2.405 \right)^2 \log(J_0).  \eqno{(69)}$$
(68) implies (61), and, if (69) holds, then (62) holds since $Z \le Z_0$ and $J \ge J_0$.  
\end{proof}

\begin{Lemma} 
Let $(a,b,c,r,s; x_1, y_1, x_2, y_2, \dots, x_N, y_N)$ and $(A,B,C,R,S; \allowbreak X_1, Y_1, \allowbreak X_2, Y_2, \allowbreak \dots, X_N, Y_N)$ be any two sets of solutions in the same family for which $(ra,sb)=(RA, SB)= 1$.  Then for every $i$ there exists a unique $j$ and for every $j$ there exists a unique $i$ such that $r a^{x_i} = R A^{X_j}$ and $s b^{y_i} = S B^{Y_j}$.  
\end{Lemma}

\begin{proof}
Suppose $(a,b,c,r,s; x_1, y_1, x_2, y_2, \dots, x_N, y_N)$ and $(A,B,C,R,S; X_1, Y_1, \allowbreak X_2, Y_2, \dots, X_N, Y_N)$ are any two sets of solutions in the same family, so that $C = kc$ for some positive rational $k$.  Suppose $(ra,sb)=(RA, SB)= 1$.  Since $(ra,sb)=1$, by the definition of family $k$ is an integer and, since $(RA, SB) = 1$, we must have $k=1$.  By Observation 3, the lemma holds.  
\end{proof}

\begin{proof}[Proof of Theorem 2:]
It suffices to prove that (59) and (60) hold when $N \ge 4$ and $(ra,sb)=1$.  

If (59) or (60) fails to hold for a given set of solutions for which $N>4$, then that set of solutions has a subset for which $N=4$ and for which (59) or (60) fails to hold.  So it suffices to consider only $N =4$.  By Lemma 16, both (59) and (60) hold for any set of solutions with $(ra,sb)=1$ in the same family as a subset (or an associate of a subset) of one of the entries listed in (39), so by Lemma 12 it suffices to consider a set of solutions $(a,b,c,r,s; x_1, y_1, x_2, y_2, x_3, y_3, x_4, y_4)$ for which one of (40), (41), or (42) holds (note that if (59) or (60) holds for the associate of a set of solutions, it holds for the set of solutions itself).  
If (42) holds, or if (40) holds with $y_1 = y_2 > 0$, we can apply Theorem 1 of \cite{Sc-St3} to the solutions $(x_2, y_2)$, $(x_3, y_3)$, and $(x_4,y_4)$ to obtain (59) and (60).  Similarly, if (41) holds with either $x_1 > 0$ or $y_2>0$, we can apply Theorem 1 of \cite{Sc-St3} to obtain (59) and (60).  So it suffices to consider only two cases: 
$$ x_1 < x_2 < x_3 < x_4, 0=y_1 = y_2 < y_3 < y_4, \eqno{(70)}$$
$$ 0 = x_1 < x_2 < x_3 < x_4, 0 = y_2 < y_1 < y_3 < y_4.  \eqno{(71)}$$

{\it Case 1: (70) holds.}

Let $(a,b,c,r,s; x_1, y_1, x_2, y_2, x_3, y_3, x_4, y_4)$ be a set of solutions satisfying (70).  Let $Z=\max(x_4, y_4)$.  Let $R= r a^{x_1}$.  

Now we can apply Lemma 13 to obtain
$$ a^{x_3-x_2} \le Z, s \le Z+1. $$
From this we have $r a^{x_3} \le r a^{x_2} Z \le (c+s) Z \le (2s+R) Z$.  Considering the solution $(x_3, y_3)$, we find $s b^{y_3} \le ra^{x_3} +c \le  (2s+R)Z + R+s$, $b^{y_3} \le (2 + R/s)Z + 1 + R/s \le (8/3) Z + 5/3$ if $R/s \le 2/3$.  Considering (2) with $(i,j)=(1,2)$, we find $R \le 2$, so that $R/s > 2/3$ implies $(R,s) = (2,1)$ or $(1,1)$. 

So assume $(R,s) \ne (2,1)$ or $(1,1)$, so that $R/s \le 2/3$.   Then, by the results in the preceding paragraph, we have 
$$a \le Z, b \le (8/3) Z +5/3, r \le R \le 2, s \le Z+1, c \le Z + 3. \eqno{(72)}$$

By (70) $\min(x_4, y_4)>0$, and, since $c \le R+s$, $(u_4, v_4) \ne (0,0)$, where $u_i$ and $v_i$ are the values of $u$ and $v$ in (1) when $(x,y) = (x_i, y_i)$.  
We can apply Lemma 15 with $(x,y)=(x_4, y_4)$ to see that (61) or (62) holds with $J =\max(a,b)$ and $d = \min(ra^{x_4}, sb^{y_4})$.  Note that $1+c/d \le \max(2,c)$.  Now, using (72) in (61) and (62), we obtain $Z < 7.4 \times 10^{14}$, and using (72) again we obtain (59) and (60), completing the proof of Theorem 2 for Case 1 when $(R,s) \ne (2,1)$ or $(1,1)$. 

Now assume $(R,s) = (2,1)$, so that we have $c =1$ or $3$.  If $c=1$, then, considering the solution $(x_2, y_2)$, we have $x_2=x_1$, a contradiction.  So $c=3$, $a^{x_2-x_1} = 2$, $x_3 > 1$, and $y_3 > 0$.  Take $3 \le i \le 4$.  We have $u_i \ne v_i$ since $c \le R+s$.  $(u_i, v_i) = (1,0)$ implies $b^{y_i} \congruent 3 \bmod 8$, while $(u_i, v_i) = (0,1)$ implies $b^{y_i} \congruent 5 \bmod 8$, so $y_i$ is odd and we must have $(u_3, v_3)= (u_4, v_4)$.  $(u_3, v_3) = (u_4, v_4) = (1,0)$ contradicts Lemma 2 of \cite{Sc-St1} since $y_4 > 1$.  So $(u_3, v_3) = (u_4, v_4)=(0,1)$, in which case, using Theorem 3 of \cite{Sc} and recalling $y_3$ and $y_4$ are both odd, we must have $(b, y_3,y_4) = (5,1,3)$ so 
$$(a,b,c,r,s; x_1, y_1, x_2, y_2, \dots, x_4, y_4) = (2,5,3,2,1;0,0,1,0,2,1,6,3)$$
or
$$(a,b,c,r,s; x_1, y_1, x_2, y_2, \dots, x_4, y_4) = (2,5,3,1,1; 1,0,2,0,3,1,7,3).$$
Both these cases certainly satisfy (59) and (60).  

Now assume $(R,s) = (1,1)$, so that $r=R=1$.  We have $c = 2$, $a^{x_2} = 3$, $x_3>1$, and $y_3 > 0$.  Take $3 \le i \le 4$.  We have $u_i \ne v_i$.  If $(u_i, v_i)=(0,1)$ and $2 \mid y_i$, then, considering the coefficients of the real term and the imaginary term of $ (1 + \sqrt{-2})^{x_i}$ modulo 9, we find that we must have $3 \mid x_i$, and, by Lemma 1 of \cite{Sc}, we find that the only possibility is $x_i=3$, giving $b=5$ and $y_i = 2$, in which case, by Theorem 7 of \cite{Sc-St2}, the only further solution is $(x,y) = (1,1)$, which contradicts (70).  So  $(u_i, v_i) = (0,1)$ implies $y_i$ odd; but then, considering the solutions $(x_3, y_3)$ and $(x_4, y_4)$ modulo 8, we find that we must have $(u_3, v_3)=(u_4,v_4)$.  If $(u_3, v_3) = (u_4, v_4) = (1,0)$ we have a contradiction to Lemma 2 of \cite{Sc-St1}.  If $(u_3, v_3)=(u_4, v_4)=(0,1)$, then, since $y_3$ and $y_4$ are both odd, we can apply Theorem 3 of \cite{Sc} to obtain a contradiction.  So Theorem 2 holds for Case 1.  

{\it Case 2: (71) holds.}

Without loss of generality we can take $a > b$.  Let $Z=\max(x_4, y_4)$.  We can apply Lemma 13 to obtain 
$$b < a \le Z, \max(r,s) \le Z+1. \eqno{(73)}$$
We have $\min(x_4, y_4) > 0$, and, since $c \le r + s b^{y_1}$, we have $(u_4, v_4) \ne (0,0)$.  Let $D = \max(r a^{x_4}, s b^{y_4} )$ and $d = \min(r a^{x_4}, s b^{y_4})$.  
Now we can apply Lemma 15 with $(x,y) = (x_4, y_4)$ to see that we must have either 
$$ Z < { \log(1+c/d) + \log(c) \over \log(2)} + 1.6901816335 \cdot 10^{10} \log(\max(r,s)) \log(a) \log(1.5 e Z), \eqno{(74a)}$$
or
$$ Z < { \log(1+c/d) + \log(c) \over \log(2)} + 22.997 \left( \log\left({Z \over \log(2)}\right) + 2.405 \right)^2 \log(a). \eqno{(74b)}$$
We first show that $\log(1 + c/d)$ is small.    
By Lemma 14,  
$$d > b^2 (c-2)/2 \ge 2 c - 4.$$
So 
$$1+{c \over d} <  1.75 \eqno{(75)}$$
when $c \ge 6$, and $1 + { c \over d} \le 1 + { c \over 8 } < 1.75$ when $c < 6$.  So (75) holds for all values of $c$. 
Now applying (73) and (75) to (74a) and (74b), and letting $K = 1.6901816335 \cdot 10^{10} $, we get either
$$ Z <  0.9 + {\log(c) \over \log(2)}  + K \log(Z+1) \log(Z) \log(1.5 e Z), \eqno{(76a)} $$
or
$$ Z <  0.9 + {\log(c) \over \log(2)} + 22.997 \left( \log\left({Z \over \log(2)}\right) + 2.405 \right)^2 \log(Z). \eqno{(76b)}$$
(76b) implies (76a), so from here on we consider only (76a).  

Suppose $c \le (Z+1)^{10^{11}}$.  
Then (76a) gives $Z < 7.9 \times 10^{14}$, so that applying (73), we get (59) and (60), completing the proof of Theorem 2 when $c \le (Z+1)^{10^{11}}$.  So from here on we can assume that, by (73),   
$$c > (Z+1)^{10^{11}} \ge (\max(r,s))^{10^{11}}. \eqno{(77)}$$
We can apply Lemma 13 to obtain 
$$x_2 a \le Z. \eqno{(78)}$$
Now $c \le r a^{x_2} + s$ so that, taking $t = {10^{11} \over 10^{11}-1 }$, we have, from (77), noting that (77) also implies $c > 10^{11} s$,
$$c < t^t a^{t x_2} $$
which gives   
$$ {\log(c) \over \log(2) } < {t \log( t) + t x_2 \log( a) \over \log(2)}. \eqno{(79)}$$
We certainly have 
$${ t \log(t) \over \log(2)} < 0.1 \eqno{(80a)}$$
and we also have 
$${t x_2 \log(a) \over \log(2)} < {x_2 a \over 1.89} \le {Z \over 1.89}  \eqno{(80b)}$$
where the first inequality follows from taking $a > 2$ and the second inequality follows from (78).  
Now we have, using (76a), (80a), and (80b), 
$$Z < 0.9 + 0.1 + { Z \over 1.89} + K \log(Z+1) \log(Z) \log(1.5 e Z)$$
from which we get  
$$0.47 Z < 1 + K \log(Z+1) \log(Z) \log(1.5 e Z), \eqno{(81)}$$
from which we obtain a bound on $Z$, which, with (73), gives (59).  

So it remains to prove (60) for Case 2 when $c > (Z+1)^{10^{11}}$.  To do this, we must significantly improve the value $0.47$ in (81), so we must significantly improve the value $1.89$ in (80b).  One way to do this is to take $a$ large: in fact, it suffices to take $a > 50$ to get the necessary improvements, and we obtain (60), completing the proof of Theorem 2 for $a > 50$.  To handle $a \le 50$, we will need two more lemmas:  

\begin{Lemma} 
Let $a>1$ and $b>1$ be relatively prime integers.  For $1 \le i \le m$, let $p_i$ be one of the $m$ distinct prime divisors of $a$.  Let $p_i^{g_i} || b^{n_i} \pm 1$, where $n_i$ is the least positive integer for which there exists a positive integer $k$ such that $\vert b^{n_i} - k p_i \vert =1$, and $\pm$ is read as the sign that maximizes $g_i$.    

Write
$$\sigma = \sum_i g_i \log(p_i)/ \log(a).$$

Then, if 
$$a^x \mid b^y \pm 1, $$
where the $\pm$ sign is independent of the above, we must have 
$$a^{x} \mid  a^{\sigma} y.$$
\end{Lemma}

\begin{proof}  Let $a = \prod_i p_i^{\alpha_i}$.  If $a^x | b^y \pm 1$, then for each $i$, $p_i^{x \alpha_i} | b^y \pm 1$, so that $p_i^{x \alpha_i -g_i} | y$ (in the case $x \alpha_i < g_i$, $p_i^{x \alpha_i - g_i}$ is a fraction that evenly divides $y$).  Thus, $y$ is divisible by 
$$\prod_i p_i^{x \alpha_i - g_i} = a^{x-\sigma}.$$
\end{proof}

\begin{Lemma}  
If $a>2$ and $(a,b) \ne (3,2)$, then, in the notation of Lemma 17, 
$$\sigma < { a \log b \over 2 \log a}.$$
If $(a,b) = (3,2)$, then $\sigma = 1$.  
\end{Lemma}

\begin{proof}  
We assume $a>2$ and $(a,b) \ne (3,2)$.  Then if $a$ is odd, $\prod_i p_i^{g_i} \le b^{\phi(a)/2} + 1 \le b^{(a-1)/2} + 1 < b^{a/2}$, verifying Lemma 18 when $a$ is odd.  If $a>4$ is even, then $\prod_i p_i^{g_i} < b^{\phi(a/2)} < b^{a/2}$ verifying the lemma in this case also.  Finally, when $a=4$, define $g$ so that $2^g || b \pm 1$, where the sign is chosen to maximize $g$.  Then the lemma holds unless ${g \log(2) \over \log(4) } \ge {4 \log(b) \over 2 \log(4)}$, that is, unless $2^g \ge b^2$, which is impossible.  

Finally, if $(a,b) = (3,2)$, $\sigma = 1$ by the definition of $\sigma$ in Lemma 17.  
\end{proof}

Returning to the proof of Theorem 2, it remains to handle Case 2 when $c > (Z+1)^{10^{11}}$ and $a \le 50$.  By (76a), (79), and (80a), 
we have 
$$Z < 1 + { t x_2 \log(a) \over \log(2)} + K \log(Z+1) \log(Z) \log(1.5 e Z). \eqno{(82)}$$
By Lemma 17 and Lemma 18 we have $a^{x_2} \le a^\sigma Z$ where $\sigma < a/2$. 
Considering the second term on the right side of (82), we have 
$$t x_2 \log(a) / \log(2) \le (t / \log(2)) (\log(Z)+ \sigma \log(a)) < 1.443 ( \log(Z) + a \log(a) /2).$$
Substituting this into (82), and recalling $a \le 50$, we get 
$$Z < 143 + 1.443 \log(Z) + K \log(Z) \log(Z+1) \log(1.5 e Z). \eqno{(83)}$$
From (83) we obtain a bound on $Z$ which, combined with (73), gives (60).      
This completes the proof of Theorem 2.  
\end{proof}

\noindent 
{\it Comment on Lemmas 17 and 18:} Lemmas 3 and 4 of \cite{Sc-St2} can be replaced by Lemmas 17 and 18 respectively, giving a shorter, simpler presentation in which Lemma 5 of \cite{Sc-St2} would be replaced by the following: 

\begin{Lemma} 
Let $a>2$, $b>1$, and $c>0$ be integers with $(a,b)=1$.  If $(-1)^u a^x + (-1)^v b^y = c$ has two solutions $(x_1, y_1, u_1, v_1)$ and $(x_2, y_2, u_2, v_2)$, with $x_2 \ge x_1 \ge 1$, $y_2 \ge y_1 \ge 1$, and $u_1$, $u_2$, $v_1$, $v_2 \in \{0,1\}$, and if further $a^{x_1} > c/2$, then 
$$x_1 < \sigma + k, $$
where $\sigma$ is defined as in Lemma 17, and $k = \frac{ 8.1 + \log \log (a) }{\log(a)}$ when $a < 5346$ and $k= 1.19408$ otherwise.  
\end{Lemma}

The proof of Lemma 19 is essentially identical to that of Lemma 5 of \cite{Sc-St2}:  we can use a result of Mignotte~\cite{Mi1} as in Proposition 4.4 of Bennett~\cite{Be}, noting that we do not need to consider the cases $(a,b)=(3,2)$ and $(a,b) = (5,2)$, since these cases are handled by the elementary methods of \cite{Pi} along with the Corollary to Theorem 2 of \cite{Sc-St1}.  
Note that in the proof of Lemma 19 we will use the bound in Lemma 18 rather than the bound in Lemma 4 of \cite{Sc-St2}: also note that we will obtain $\frac{ y_2 \log(b) } { \log(c)} > 10.519$ rather than $\frac{ y_2 \log(b) }{\log(c)} > 34$ as in \cite{Sc-St2}.  Lemmas 3, 4, and 5 of \cite{Sc-St2} are used to prove Theorem 3 of \cite{Sc-St2}; if we use instead Lemmas 17, 18, and 19 of the present paper, both the formulation and the proof of Theorem 3 of \cite{Sc-St2} remain completely unchanged.

\section{The Case $N \le 3$ with $\gcd(ra,sb)=1$  }  

The question remains: what can be said of cases for which $N \le 3$?  

It follows directly from Theorem 2 of \cite{Sc-St3} that cases of exactly two solutions to (1) with $(ra,sb)=1$ are commonplace and easy to construct, even if we restrict consideration to basic forms only (here we are allowing $N=2$, and noting that we can easily adjust the two solutions to get $\min(x_1, x_2) = \min(y_1, y_2) = 0$).    

When $N=3$, we find many examples.  Here we list several types of sets of solutions, each one of which generates an infinite number of families giving three solutions to (1).  We list these sets of solutions in the form $(a,b,c,r,s;x_1, y_1, x_2, y_2, x_3, y_3)$:  

$$(a, {a^{kd} + (-1)^{u+v} \over a^d + (-1)^{u}}, { a^d b + (-1)^{u+v+1} \over h}, {b+(-1)^v \over h}, {a^d+(-1)^{u} \over h}; 0, 1, d, 0, kd, 2) \eqno{(84)}$$
where $a$ and $b = {a^{kd} + (-1)^{u+v} \over a^d + (-1)^{u}}$ are integers greater than 1, $d$ and $k$ are positive integers, $h=\gcd(a^d+(-1)^u, b+(-1)^v)$, and $u$ and $v$ are in the set $\{ 0,1 \}$.  When $u=0$, we take $k-v$ odd; when $(u,v) = (1,1)$, we take $a^d \le 3$.  When $a=d=2$ and $(u,v)=(1,1)$, we can take $k$ to be a half integer.  When $k=2$ and $u-v$ is odd, the same choice of $(a,b,r,s)$ as in (84) gives the additional set of solutions 
$$(a, a^d +(-1)^v, {2 a^d + (-1)^v \over h} ,{a^d +(-1)^v 2 \over h}, {a^d + (-1)^{v+1} \over h}; 0,0,d,1,3d,3).  \eqno{(85)}$$
Other sets of solutions can be constructed with specified values of $a$.  For example, when $a=3$ we have 
$$(3, {3^g +(-1)^v \over 2}, {3^{g+1} + (-1)^{v}  \over 2^{2+v-\alpha} }, { 3(3^{g-1}+(-1)^{v}) \over 2^{2+v-\alpha}}, 2^{1-v+\alpha}; 0,1, 1,0,2g,3) \eqno{(86)}$$
where $v \in \{0,1\}$, $g$ is a positive integer, $\alpha = 0$ when $2 \mid g-v$, $\alpha=1$ when $g$ is odd and $v=0$, and $\alpha = 2$ when $g$ is even and $v=1$.  Note that, when $g$ is odd and $v=1$, (86) corresponds to (7) in Theorem 1 of \cite{Sc}.  Note also that the cases $(g,v) = (1,0)$ and $(2,1)$ correspond to the cases $(3,2,5)$ and $(3,2,13)$ in the exceptional cases of Theorem 7 of \cite{Sc-St1}.  

When $a=2$, we have 
$$(2, 2^g + (-1)^v, 2^{g} +(-1)^{v+1}, 2, 1; 0,1,g-1,0,g,1) \eqno{(87)}$$
where $v \in \{0,1\}$ and $g$ is a positive integer.  Note that (87) corresponds to the cases $(M, M+2)$ and $(F, F-2)$ at the end of the formulation of Theorem 7 of \cite{Sc-St1}.  

Also, it is easy to construct sets of solutions for which $x_1 = y_1 = y_2 = 0$.  For example, we have, for $a$ even and $x>0$, 
$$(a, 2 a^x \pm 1, a^x \pm 1, 2, a^x \mp 1; 0,0,x,0,2x,1), \eqno{(88)}$$
or, more generally,
$$(a, {2a^{x_3} +(-1)^{t+w+1}a^{x_2} + (-1)^{w+1} \over  a^{x_2} + (-1)^{t+1} } , { a^{x_2} +(-1)^{t} \over 2^{m} } , 2^{1-m}, { a^{x_2} +(-1)^{t+1} \over 2^m }; 0,0,x_2, 0, x_3, 1), \eqno{(89)}$$
where $x_2 >0$ and $x_3>0$ are chosen so that $a^{x_3} \congruent (-1)^w \bmod {a^{x_2} +(-1)^{t+1} \over 2^m}$, $t \in \{ 0,1 \}$, $w \in \{0,1\}$, and $m=1$ or 0 according as $a$ is odd or even.    

Since each of (84--89) generates an infinite number of families for which $N=3$, we have, after Theorem 2, 

\begin{Theorem} 
There are an infinite number of cases of exactly three solutions to (1) with $(ra, sb)=1$, even if we consider only sets of solutions in basic form. 
\end{Theorem}
 
\noindent 
{Comment on anomalous cases (not necessarily with $(ra,sb)=1$):} 

If we exclude from consideration any set of solutions in the same family as a set of solutions given by either (57) or (58) or any of (84--89), then we are aware of only 14 essentially different cases of $(a,b,c,r,s)$ giving exactly three solutions to (1), the largest of which is 
$$(56744, 1477,  83810889,  1478, 56743;  0,1,1,0, 3,4).  \eqno{(90)}$$
Easily derived from (90) is $ (56745,  1477,  41906182,  739, 28373;  0,1,  1,0,  3,4)$.


\begin{thebibliography}{1}

\bibitem{Be}
M. Bennett, 
\newblock On some exponential equations of S. S. Pillai, 
\newblock {\it Canadian Journal of Mathematics}, {\bf 53}, no. 5, (2001), 897--922.

\bibitem{BBM}
Y.~F. Bilu, Y. Bugeaud, M. Mignotte, 
\newblock {\it Catalan's Equation}, book in preparation.  

\bibitem{BL}
Y. Bugeaud, F. Luca, 
\newblock On Pillai's Diophantine equation, 
\newblock {\it New York J. of Math.}, {\bf 12}, (2006) 193--217. (electronic)

\bibitem{GLS} 
R.~K. Guy, C.~B. Lacampagne, J.~L. Selfridge, 
\newblock Primes at a glance, 
\newblock {\it Math. Comp.}, {\bf 48}, (1987), 183--202. 

\bibitem{HT}
B. He, A. Togb\'e, 
\newblock On the number of solutions of the exponential Diophantine equation $a x^m - b y^n = c$,
\newblock{\it Bull. Aust. Math. Soc.}, {\bf 81},  (2010), 177--185.

\bibitem{Le}
M. Le,  
\newblock A note on the diophantine equation $a x^m - b y^n = k $, 
\newblock {\it Indag. Math.} (N. S.) {\bf 3},  June 1992, 185--191.

\bibitem{Ma}
E.~M. Matveev,  
\newblock An explicit lower bound for a homogeneous rational linear form in logarithms of algebraic numbers. II. , 
\newblock {\it Izv. Math.} {\bf 64} (2000), 1217--1269. 

\bibitem{Mi1}
M. Mignotte,
\newblock A corollary to a theorem of Laurent-Mignotte-Nesterenko, 
\newblock {\it Acta Arith.}, {\bf 86} (1998), 101--111.

\bibitem{Mi2}
M. Mignotte, 
\newblock A kit on linear forms in three logarithms, in preparation, February 7, 2008.   {\textnormal http://www-irma.u-strasbg.fr/\~ bugeaud/travaux/kit.pdf}

\bibitem{Pi}
S. S.  Pillai, 
\newblock On the equation $2^x-3^y = 2^X + 3^Y$,
\newblock {\it Bull. Calcutta Soc.} , {\bf 37}, (1945), 15--20.

\bibitem{Sc}
R. Scott, 
\newblock On the Equations $p^x-b^y = c$ and $a^x+b^y=c^z$, 
\newblock {\it Journal of Number Theory}, {\bf 44}, no. 2 (1993), 153--165. 

\bibitem{Sc-St1}
R. Scott, R. Styer, 
\newblock On $p^x - q^y = c$ and related three term exponential Diophantine equations with prime bases, 
\newblock {\it Journal of Number Theory}, {\bf 105} no. 2 (2004), 212--234.

\bibitem{Sc-St2}
R. Scott, R. Styer, 
\newblock On the generalized Pillai equation $\pm a^x \pm b^y = c$, 
\newblock {\it Journal of Number Theory}, {\bf 118} (2006), 236--265.

\bibitem{Sc-St3}
R. Scott, R. Styer, 
\newblock The generalized Pillai equation $\pm r a^x \pm s b^y = c$, 
\newblock accepted by {\it Journal of Number Theory}.

\bibitem{Sc-St5}
R. Scott, R. Styer, 
\newblock Handling a large bound for a problem on the generalized Pillai equation $\pm r a^x \pm s b^y = c$,
\newblock preprint.

\bibitem{Sh}
T.~N. Shorey,
\newblock On the equation $a x^m - b y^n = k$
\newblock  {\it  Nederl. Akad. Wetensch. Indag. Math.} {\bf  48 } (1986), no. 3, 353--358. 

\bibitem{St1}
R. Styer, 
\newblock Small two-variable exponential Diophantine equations, 
\newblock {\it Mathematics of Computation}, {\bf 60}, No. 202, 1993, 811--816.

\bibitem{W}
M. Waldschmidt, 
\newblock {Perfect powers: Pillai's works and their developments, Arxiv preprint arXiv:0908.4031,
August 27, 2009.}

\end{thebibliography}
\end{document}